\documentclass[11pt]{article}
\usepackage{epigamath}

\setpapertype{A4}

\usepackage[english]{babel}

\usepackage[dvips]{graphicx}     

\title{Correspondences between convex geometry and complex geometry}
\titlemark{Correspondences between convex geometry and complex geometry}
\author{Brian Lehmann and Jian Xiao}
\authoraddresses{
\authordata{Brian Lehmann}{\firstname{Brian} \lastname{Lehmann}\\
\institution{Department of Mathematics, Boston College, Chestnut Hill, MA 02467, USA}\\
\email{lehmannb@bc.edu}}

\authordata{Jian Xiao}{\firstname{Jian} \lastname{Xiao}\\
\institution{Department of Mathematics, Northwestern University, Evanston, IL 60208, USA}\\
\email{jianxiao@math.northwestern.edu}}
}
\authormark{B. Lehmann and J. Xiao}
\date{\vspace{-5ex}} 
\journal{\'Epijournal de G\'eom\'etrie Alg\'ebrique} 
\acceptation{Received by the Editors on September 20, 2016, and in final form on March 15, 2017. \\ Accepted on June 19, 2017.}

\newcommand{\articlereference}{Volume 1 (2017), Article Nr. 6}

\newtheorem{thrm}{Theorem}[section]
\newtheorem{lem}[thrm]{Lemma}
\newtheorem{cor}[thrm]{Corollary}
\newtheorem{prop}[thrm]{Proposition}

\newtheorem{defn}[thrm]{Definition}

\newtheorem{exmple}[thrm]{Example}
\newtheorem{rmk}[thrm]{Remark}
\newtheorem{ques}[thrm]{Question}

\DeclareMathOperator{\Eff}{\overline{Eff}}
\DeclareMathOperator{\Nef}{Nef}
\DeclareMathOperator{\Mov}{Mov}

\DeclareMathOperator{\vol}{vol}

\DeclareMathOperator{\Sym}{Sym}

\DeclareMathOperator{\codim}{codim}

\DeclareMathOperator{\Image}{Image}
\DeclareMathOperator{\PSH}{PSH}

\begin{document}


\maketitle



\begin{prelims}

\vspace{-0.1cm}

\def\abstractname{Abstract}
\abstract{We present several analogies between convex geometry and the theory of holomorphic line bundles on smooth projective varieties or K\"ahler manifolds. We study the relation between positive products and mixed volumes. We define and study a Blaschke addition for divisor classes and mixed divisor classes, and prove new geometric inequalities for divisor classes. We also reinterpret several classical convex geometry results in the context of algebraic geometry: the Alexandrov body construction is the convex geometry version of divisorial Zariski decomposition; Minkowski's existence theorem is the convex geometry version of the duality between the pseudo-effective cone of divisors and the movable cone of curves.}



\vspace{0.05cm}

\languagesection{Fran\c{c}ais}{%

\textbf{Titre. Correspondances entre g\'eom\'etrie convexe et g\'eom\'etrie complexe}
\commentskip
\textbf{R\'esum\'e.}
Nous pr\'esentons plusieurs analogies entre la g\'eom\'etrie convexe et la th\'eorie des fibr\'es en droites sur les vari\'et\'es projectives lisses ou k\"ahl\'eriennes. Nous \'etudions les relations entre produits positifs et volumes mixtes. Nous d\'efinissons et \'etudions une somme de Blaschke pour les classes de diviseurs (mixtes ou non), et nous \'etablissons de nouvelles in\'egalit\'es g\'eom\'etriques pour les classes de diviseurs. Nous r\'einterpr\'etons \'egalement plusieurs r\'esultats classiques de g\'eom\'etrie convexe en termes de g\'eom\'etrie alg\'ebrique : la construction du corps d'Alexandrov est la version g\'eom\'etrie convexe de la d\'ecomposition de Zariski divisorielle; de m\^eme, le th\'eor\`eme d'existence de Minkowski est le pendant en g\'eom\'etrie convexe de la dualit\'e entre c\^one des diviseurs pseudo-effectifs et celui des courbes mobiles.}

\end{prelims}



\newpage

\setcounter{tocdepth}{1}
\tableofcontents
\vspace*{2cm}

\section{Introduction}

Many deep results concerning holomorphic line bundles in complex geometry are inspired by results in convex geometry.  For example, the Khovanskii-Teissier inequalities for $(1,1)$-classes are inspired by the Alexandrov-Fenchel
inequalities for convex bodies.  Such results were first explored in \cite{teissier82} and \cite{khovanskii88};   \cite{gromov1990convex} also contains many nice discussions on these relations.

In this paper we extend the dictionary between the geometry of line bundles in complex geometry and the theory of convex bodies.  First, we prove new results about divisors motivated by results in convex geometry.  These results strengthen and clarify existing geometric inequalities in the literature.  A key role is played by the refined structure of the movable cone of curves we obtained in our previous work \cite{lehmann2016positiivty}.  Second, we reinterpret several classical convex geometry results in the context of algebraic geometry.  Although these convex geometry results are well understood, we feel that complex geometers could benefit from this discussion.

The outline of this correspondence is as follows.  On the complex geometry side, let $X$ be a smooth projective variety (or compact K\"ahler manifold) of dimension $n$.  We have the vector space $H^{1,1}(X,\mathbb{R})$ whose integral points parameterize Chern classes of holomorphic line bundles.  Inside of this space, we have cones $\Eff^{1}(X)$, $\Mov^{1}(X)$ which represent $(1,1)$-classes satisfying various positivity conditions.  On the convex geometry side, we consider $$V:=C^0 (\mathbb{S}^{n-1})/L(\mathbb{S}^{n-1})$$
where $\mathbb{S}^{n-1}$ is the $(n-1)$-dimensional sphere, $C^0 (\mathbb{S}^{n-1})$ is the space of continuous functions and $L(\mathbb{S}^{n-1})$ is the space of linear functions. We consider two positive cones in $V$: $\mathcal{C}_p$ is the cone generated by positive continuous functions and $\mathcal{C}$ is the cone generated by positive convex and homogeneous of degree one functions. We then study the following correspondences:\bigskip

\begin{center}
\begin{tabular}{ c | c | c  }
   & Complex Geometry & Convex Geometry \\ \hline
   Section \hyperref[section positiveproduct]{3} &  \begin{tabular}{@{}c@{}} movable classes \\  and the positive products \end{tabular}  &  $\mathcal{C}$ and the mixed volumes \\ \hline
  Section \hyperref[section alex body]{4} & \begin{tabular}{@{}c@{}} pseudo-effective classes and \\ the divisorial Zariski decomposition \end{tabular}  &  \begin{tabular}{@{}c@{}} $\mathcal{C}_{p}$ and\\ the Alexandrov body construction\end{tabular}   \\ \hline
  Section \hyperref[morseineq section]{5} & Morse inequality for $(\Eff^1 (X), \vol)$ & Morse inequality for $(\mathcal{C}_p, \vol)$  \\ \hline
  Section \hyperref[section minkowski]{6} & \begin{tabular}{@{}c@{}} movable curve classes \\  and duality of cones \end{tabular} & \begin{tabular}{@{}c@{}} mixed area measures \\  and Minkowski's existence theorem \end{tabular}  \\ \hline
  Section \hyperref[section blachke]{7} & addition of movable curve classes & \begin{tabular}{@{}c@{}} Blaschke addition and \\   the Kneser-S\"uss inequality  \end{tabular} \\  \hline
  Section \hyperref[section loomis]{8} & fibration and restricted volumes & projection and projection volumes  \\
\end{tabular}
\end{center}
\bigskip

In this dictionary, convex bodies with non-empty interior correspond to big movable classes. We can phrase several classical convex geometry results due to Minkowski and Alexandrov in terms of algebraic geometry; for example, the divisorial Zariski decomposition is a variant of the Alexandrov body construction.

We highlight a couple particular results.  The following proposition\footnote{The first part of this proposition has been implicitly contained in \cite{fl14} and \cite{lehmann2016positiivty}, and was also mentioned to the second author by S.~Boucksom.} corresponds to the linearity of the $(n-1)$st mixed volume of polytopes.

\begin{prop}
Let $X$ be a smooth projective variety of dimension $n$. Assume that $L$ is pseudo-effective and $M$ is movable, then
\begin{equation*}
  \langle L^{n-1}\rangle \cdot M = \langle L^{n-1}\cdot M\rangle.
\end{equation*}
However, there exist projective manifolds such that $\langle L_1 \cdot L_2 \cdot...\cdot L_{n-1}\rangle \cdot M \neq \langle L_1 \cdot L_2 \cdot...\cdot L_{n-1}\cdot M\rangle$ if the pseudo-effective classes $L_i$ are different.
\end{prop}

We also discuss Brascamp-Lieb type inequalities on smooth projective varieties; for example, we prove the following Loomis-Whitney inequality for $(1,1)$ classes.

\begin{thrm}[Theorem \ref{thrm loomis}]
Let $X$ be a smooth projective variety of dimension $n$.  Suppose that $\{ f_{i}: X \dashrightarrow C_{i} \}_{i=1}^{n}$ is a collection of dominant rational maps onto curves, and assume that the induced map $X \dashrightarrow \prod_{i=1} ^n C_{i}$ is birational. Then for any big $(1,1)$ class $L$, we have
\begin{equation*}
\left( \vol(L)/n! \right)^{n-1}  \leq \prod_{i=1}^{n} \left( \vol_{f_i} (L)/(n-1)! \right),
\end{equation*}
where $\vol_{f_i} (L)$ is the restricted volume on a general fiber of $f_i$ (up to a resolution). Moreover, if all the $f_{i}$ are morphisms, then we obtain equality if and only if the positive part of $L$ is numerically proportional to a sum of fibers of the $f_{i}$.
\end{thrm}

While our discussion of the fields of complex and convex geometry will be purely by analogy, one can sometimes make the relationship between line bundles and polytopes quite explicit.  There is a well-understood way to construct a polytope (up to translation) to a line bundle in the toric setting.  More generally, one has an analogous construction for arbitrary varieties known as the Okounkov or Newton-Okounkov body (see e.g. \cite{lm09}, \cite{kk12}).  Unfortunately, certain geometric features of Okounkov bodies rely heavily on the choice of flag.   Thus it is not clear whether the dictionary we present here can be enriched to the polytope setting while still capturing the full strength of our results.  In particular, it is interesting to ask the following question:

\begin{ques} \rm
Let $X$ be a compact K\"ahler manifold (or smooth projective variety) of dimension $n$.  After possibly fixing some data on $X$, is there a ``canonical'' injective linear map
\begin{equation*}
  i: H^{1,1} (X, \mathbb{R}) \rightarrow V:=C^0 (\mathbb{S}^{n-1})/L(\mathbb{S}^{n-1}),
\end{equation*}
or equivalently, surjective linear map
\begin{equation*}
  i^*: V^* \rightarrow H^{n-1,n-1} (X, \mathbb{R}),
\end{equation*}
such that the positive cones of $(1,1)$ classes are canonically embedded in the positive cones of functions, the map $i$ preserves positivity and the pairing between $H^{1,1} (X, \mathbb{R})$ and $H^{n-1,n-1} (X, \mathbb{R})$ is induced by the pairing between $V$ and $V^*$ via the maps $i, i^*$?
\end{ques}

\subsection*{Acknowledgement}
We would like to thank J.~P.~Demailly, M.~Jonsson and V.~Tosatti for several comments and suggestions on this work, and R.~Lazarsfeld, M.~Musta{\c t}{\u a}, and B. Totaro for pointing out an error in an earlier version.   The first author was supported by an NSA Young Investigator Grant.  The second author was supported in part by CSC fellowship and Institut Fourier.

\section{Preliminaries}\label{section preliminaries}
\subsection{Positive classes}

Let $X$ be a compact K\"ahler manifold of dimension $n$.  We will let $H^{1,1}(X, \mathbb{R})$ denote the real Bott-Chern cohomology group of bidegree $(1,1)$. Note that the Bott-Chern cohomology coincides with de Rham cohomology on a compact K\"ahler manifold, due to the $\partial\bar \partial$-lemma. A Bott-Chern class is called pseudo-effective if it contains a $d$-closed positive current.  A pseudo-effective class $L$ is called movable if for any irreducible divisor $Y$ the Lelong number  $\nu(L,y)$ (or minimal multiplicity as in \cite{Bou04}) vanishes at a very general point $y \in Y$.  We will be interested in the following cones in $H^{1,1}(X,\mathbb{R})$:
\begin{itemize}
\item $\Eff^{1}(X)$: the cone of pseudo-effective $(1,1)$-classes;
\item $\Mov^{1}(X)$: the cone of movable $(1,1)$-classes.
\end{itemize}
Classes in the interior of $\Eff^{1}(X)$ are known as big classes; these are exactly the classes admitting a K\"ahler current (i.e. strictly positive current), or equivalently, having strictly positive volumes. As a consequence of Demailly's regularization theorem \cite{demailly1992regularization}, any big class contains a K\"ahler current with analytic singularities.
Let $H^{n-1,n-1}(X, \mathbb{R})$ denote the real Bott-Chern cohomology group of bidegree $(n-1,n-1)$. We will be interested in the following cone in $H^{n-1,n-1}(X,\mathbb{R})$:
\begin{itemize}
\item $\Mov_{1}(X)$: the cone of movable $(n-1,n-1)$-classes.
\end{itemize}
Recall that $\Mov_1 (X)$ is the closed cone generated by classes of the form $\mu_* (\widetilde{A}_1 \cdot ... \cdot \widetilde{A}_{n-1})$, where $\mu$ is a modification and the $\widetilde{A}_i$ are K\"ahler classes upstairs, see e.g. \cite{BDPP13}. An irreducible curve $C$ is called movable if it is a member of an algebraic family that covers $X$.

\begin{rmk} \rm
Let $\omega$ be a fixed K\"ahler class on $X$ and let $L$ be a $(1,1)$ class.  By \cite{Bou04} the condition $L \in \Mov^1 (X)$ is equivalent to requiring that for any $\delta>0$ there is a modification $\mu: \widehat{X}\rightarrow X$ such that $L+\delta \omega = \mu_* \widehat{\omega}_\delta$, where $\widehat{\omega}_\delta$ is a K\"ahler class on $\widehat{X}$. In particular, when $n=2$, we have $\Mov^1 (X) = \Mov_1 (X)$.
\end{rmk}

The most important conjecture concerning these cones is Demailly's conjecture on the (weak) transcendental Morse inequality.  This conjecture was recently proved for projective manifolds in \cite{nystrom2016duality}.  When working with positive cones, it will at times be useful to restrict ourselves to this setting. In particular, by \cite{BDPP13, nystrom2016duality}, for projective manifolds we have
\begin{equation*}
  \Eff^1 (X)^* =\Mov_1 (X).
\end{equation*}

\begin{rmk} \rm
By \cite{Cut13} and \cite[Section 2.2.2]{fl14}, the results of \cite{BFJ09} can be extended to an arbitrary algebraically closed field. Thus it is also possible to work in the purely algebraic setting, i.e.~with smooth algebraic varieties over an arbitrary algebraically closed field
(see also \cite[Section 2.3]{lehmxiao16convexity}).
In this case one should replace $H^{1,1}(X,\mathbb{R})$ by the space $N^{1}(X)$ of $\mathbb{R}$-Cartier divisors up to numerical equivalence.  All our results hold equally well with this change.
\end{rmk}

\begin{rmk} \rm
As already noted in our previous works \cite{lehmxiao16convexity, lehmann2016positiivty}, the main results of these papers extend to the transcendental situation whenever we have the (weak) transcendental Morse inequality. In particular, the transcendental versions hold on projective manifolds.
\end{rmk}

\subsection{Positive functions on the sphere}
Let $\mathfrak{K}^n$ denote the set of convex bodies (compact convex subsets) with non-empty interior in $\mathbb{R}^n$. For two convex bodies $K, L$, we have the Minkowski sum $K+L:=\{x+y|\ x\in K, y\in L\}$ and it is clear $K+L \in \mathfrak{K}^n$. In this way, $(\mathfrak{K}^n, +)$ is a commutative semi-group. For $\lambda\geq 0$,  the scaling operation is $\lambda K:=\{\lambda x|\ x\in K\}$.

Associated with a convex body $K$ is its support function, $h(K, \cdot)$, defined on $\mathbb{R}^n$ by
\begin{align*}
h(K, x):=\max\{x\cdot y|y\in K\}.
\end{align*}
We shall write $h_K (\cdot)$, rather than $h(K, \cdot)$, for the support function of $K$, and we will usually consider only the restriction of the support function to $\mathbb{S}^{n-1}$. So $h_K$ is a continuous function on $\mathbb{S}^{n-1}$.
Let $C^0 (\mathbb{S}^{n-1})$ denote the real valued continuous function space on the sphere. Then the support functions induce a map:
\begin{align*}
\Phi: &\mathfrak{K}^n \rightarrow C^0 (\mathbb{S}^{n-1})\\
&K \mapsto h_K (\cdot).
\end{align*}

Indeed, by the embedding theorem for spaces of convex sets (see e.g. \cite{hormanderembedding}, \cite[Chapter 1]{schneiderconvex}), $\Phi$ is an isometrical isomorphism onto its image if we endow the spaces $\mathfrak{K}^n$ and $C^0 (\mathbb{S}^{n-1})$ with suitable metrics. Let $d_H$ denote the Hausdorff metric on $\mathfrak{K}^n$, that is,
\begin{align*}
d_H (K, L)=\min\{\epsilon>0| K\subseteq L+\epsilon B\ \&\ L\subseteq K+\epsilon B\},
\end{align*}
where $B$ is the unit $n$-ball. We endow $C^0 (\mathbb{S}^{n-1})$ with maximal norm $|\cdot|_{\infty}$, so
\begin{align*}
d_H (K, L)=|h_K -h_L|_{\infty}.
\end{align*}
Thus $\Phi: (\mathfrak{K}^n, d_H) \rightarrow (\Phi(\mathfrak{K}^n), |\cdot|_{\infty})$ is an isometrical isomorphism. More importantly,
by the above embedding $\Phi$, we can identify $\mathfrak{K}^n$ with its image $\Phi (\mathfrak{K}^n)\subset C^0 (\mathbb{S}^{n-1})$.  Thus we will freely pass back and forth between convex bodies and their support functions on the sphere.

\subsubsection{Quotient space}
By the embedding theorem for convex bodies,  $\mathfrak{K}^n$ can be seen as a convex cone in an infinite dimensional Banach space. For $a\in \mathbb{R}^n$ it is clear that $h_{K+a}= h_K + l_a $ where $l_a (x) = a \cdot x$ is a linear function. Comparing the equality characterization of Brunn-Minkowski inequalities in convex geometry with the one in algebraic geometry (see Section \hyperref[section bm ineq]{2.C}), it is convenient to identify two convex bodies if they differ by a translation. This leads to the quotient space of $C^0 (\mathbb{S}^{n-1})$ by linear functions. We denote $L(\mathbb{S}^{n-1})$ to be the vector space whose elements are the restriction of linear functions to $\mathbb{S}^{n-1}$. It is obvious that $\dim L(\mathbb{S}^{n-1}) =n$. We consider the quotient space
\begin{equation*}
  V:= C^0 (\mathbb{S}^{n-1})/L(\mathbb{S}^{n-1}).
\end{equation*}
We endow $V$ with the quotient norm, that is, for $f\in V$ we define $$||f||=\inf_{l\in L(\mathbb{S}^{n-1})} |f-l|_{\infty}.$$
It is a standard fact in functional analysis that $(V, ||\cdot||)$ is also an infinite dimensional Banach space.

\begin{rmk} \rm
Unlike the finite dimensional situation, the Banach space $V$ is not reflexive. If we identify $V$ with its canonical image in $V^{**}$, then $V\subsetneq V^{**}$. This means that not every continuous linear functional on $V^*$ is given by a continuous function.
\end{rmk}

\subsubsection{Positive cones on $V$}
We will denote by $\mathcal{C}$ the image of $\mathfrak{K}^n$ in $V$, so that $\mathcal{C}$ is a convex cone in $V$.  We will also let $\mathcal{C}_{p} \subset V$ denote the convex cone generated by non negative continuous functions.  Note that $\mathcal{C} \subset \mathcal{C}_{p}$.  Motivated by the toric case, we are guided by the following analogies:
\begin{equation*}
\Eff^{1}(X) \leftrightsquigarrow \mathcal{C}_{p}, \qquad \qquad \Mov^{1}(X) \leftrightsquigarrow \mathcal{C}.
\end{equation*}

\subsection{Volume and Brunn-Minkowski inequalities}\label{section bm ineq}
In convex geometry, we have the following fundamental geometric inequality (see e.g. \cite[Chapter 7]{schneiderconvex}).

\begin{thrm}
Let $K, L \subseteq \mathbb{R}^n$ be two convex bodies with non-empty interior, then
\begin{equation*}
  \vol(K+L)^{1/n} \geq \vol(K)^{1/n} +\vol(L)^{1/n}
\end{equation*}
with equality if and only if $K, L$ are homothetic.
\end{thrm}

In the complex geometry setting, the volume of a class $L \in \Eff^{1}(X)$ is defined to be
\begin{equation*}
\vol(L) = \sup_{T} \int_{X} T_{ac}^{n}
\end{equation*}
for $T$ ranging over the closed positive $(1,1)$-currents in $L$, where $T_{ac}$ denotes the absolutely continuous part (as in \cite{bouck02linebdl}).  In particular, for divisor classes, this analytic definition coincides with the usual one.

By extending a previous result of \cite{BFJ09} from big nef classes to big movable classes,  we proved in our previous work \cite{lehmann2016positiivty}:

\begin{thrm}\label{lx brunnmink}
Let $X$ be a smooth projective variety of dimension $n$. Assume that $L, M$ are two big and movable $(1,1)$-classes, then
\begin{equation*}
 \vol(L+M)^{1/n} \geq \vol(L)^{1/n} +\vol(M)^{1/n}
\end{equation*}
with equality if and if $L, M$ are proportional.
\end{thrm}

\begin{rmk} \rm
For non-projective K\"ahler manifolds, by \cite{FX14} we also have a Brunn-Minkowski inequality and its equality characterization for big and nef $(1,1)$ classes.
\end{rmk}

\section{Positive products and mixed volumes}\label{section positiveproduct}

For toric varieties, a big movable divisor corresponds to a convex polytope (up to translation).  The most important structure on the set of polytopes is the mixed volume.  \cite{BFJ09} shows how, in the toric setting, the mixed volume corresponds to the positive product construction defined below (up to a constant factor).

\subsection{Mixed volumes and area measures}

In the convex geometry setting, let $K_1, ..., K_r \in \mathfrak{K}^n$ be convex bodies.  Then there is a polynomial relation
\begin{equation*}
  \vol(t_{1} K_1  +...+ t_{r} K_r) =\sum_{i_1 + ... + i_r =n} \frac{n!}{i_1 !i_2 !...i_r !} V(K_1 ^{i_1}, ...,K_r ^{i_r}) t_{1}^{i_{1}} \ldots t_{r}^{i_{r}},
\end{equation*}
where $t_i \in \mathbb{R}_+$ and the $V(K_1 ^{i_1}, ...,K_r ^{i_r})$ are coefficients defining the mixed volumes. In particular, if $K_1 =K_2=...=K_r =K$ (up to translations), then the mixed volume is just $\vol(K)$.

There is a variant of this construction that will frequently be useful for us.  The area measure $S(K,...,K;\cdot)$ of a convex body $K$ is, roughly speaking, the pushforward of the $(n-1)$ dimensional Hausdorff measure on $\partial K$ by its Gauss map (see \cite[Chapter 4]{schneiderconvex} for the precise definition). Alternatively, the (mixed) area measures can be defined through mixed volumes and Riesz representation theorem in a very simple way. We follow the elegant discussion in \cite{lutwak1988intersection}. Suppose that $K_1,..., K_{n-1}$ are fixed convex bodies with non empty interior.  Via the embedding map $\Phi$, the function
\begin{equation*}
  l(h_L):= V(K_1, ..., K_{n-1}, L)
\end{equation*}
can be extended in a unique way to a continuous functional on $C^0 (\mathbb{S}^{n-1})$. Here, $V(K_1, ..., K_{n-1}, L)$ is the mixed volume of convex bodies (see \cite[Chapter 5]{schneiderconvex}). Thus there exists a Borel measure, denoted by $S(K_1,..., K_{n-1}; \cdot)$, such that
\begin{align*}
l(h_L)=V(K_1,...,K_{n-1},L)
=\int_{\mathbb{S}^{n-1}}h_L(x)dS(K_1,...,K_{n-1};x).
\end{align*}
The mixed area measure $S(K_1,..., K_{n-1}; \cdot)$ has the origin as its centroid and does not concentrate on any great subsphere.

\subsection{Positive products}

Let $X$ be a compact K\"ahler manifold of dimension $n$. Before going on, we first recall the definition of positive currents with minimal singularities. Let $\{\theta\}$ be a pseudo-effective $(1,1)$ class, where $\theta$ is a smooth representative. Then a positive current $T=\theta + i\partial \bar \partial \varphi$ has minimal singularities if for any $u \in \PSH(X, \theta)$ there exists some constant $c$ such that $u\leq \varphi + c$.  Note that such a current $T_{\min}$ always exists: let
$$V_{\theta} := \sup \{\varphi \in \PSH(X, \theta)|\ \varphi \leq 0\},$$
then $\theta+ dd^c V_{\theta}$ is a positive current with minimal singularities. There may be many positive currents with minimal singularities in a class.

Assume that $L_1, ..., L_r \in H^{1,1} (X, \mathbb{R})$ are big classes, i.e. every class $L_i$ contains a K\"ahler current. By the theory developed in \cite{BEGZ10MAbig}, we can associate to $L_1, ..., L_r$ a positive class in $H^{r,r} (X, \mathbb{R})$, denoted by $\langle L_1 \cdot ...\cdot L_r\rangle$. It is defined as the class of the non-pluripolar product of positive current with minimal singularities, that is,
\begin{equation*}
  \langle L_1 \cdot ...\cdot L_r\rangle := \{\langle T_{1, \min} \wedge ...\wedge T_{r, \min}\rangle\}
\end{equation*}
where $\langle T_{1, \min} \wedge ...\wedge T_{r, \min}\rangle$ is the non-pluripolar product.  It is proved in \cite{BEGZ10MAbig} that the positive product $\langle L_1 \cdot ...\cdot L_r\rangle$ does not depend on the choices of currents with minimal singularities. Moreover, if the $L_i$ are nef, then $\langle L_1 \cdot ...\cdot L_r\rangle=L_1 \cdot ...\cdot L_r$ is the usual cohomology intersection.

\begin{rmk} \rm
As proved in \cite[Proposition 1.10]{principato2013mobile}, it is also not hard to see that the construction in \cite[Theorem 3.5]{BDPP13} and the construction in \cite{BEGZ10MAbig} are the same.
This follows directly from the maximal integration property of $\langle L_1,..., L_r\rangle$ against smooth positive closed $(n-r, n-r)$-forms (see \cite[Proposition 1.18 and Proposition 1.20]{BEGZ10MAbig}).
\end{rmk}

\begin{rmk} \label{algversusana} \rm
Let $X$ be a smooth projective variety of dimension $n$, and assume that $L_1, ..., L_r$ are real big divisor classes. Then we can also define $\langle L_1,..., L_r\rangle$ by the algebraic construction in \cite{BFJ09}.  Since the positive product constructions in \cite{BFJ09} are defined by taking limits in the space $N_{n-r} (X)$ of $\mathbb{R}$-cycles up to numerical equivalence (not in $H^{r,r}(X, \mathbb{R})$), it might be different from the analytic construction, depending on the validity of the Hodge Conjectures.  However, the resulting positive products are the same for $r=1$ and $r=n-1$ because of the Lefschetz theorem on $(1,1)$ and $(n-1, n-1)$ classes. And the equality is obvious for $r=n$.
\end{rmk}

In this section, we discuss the analogue between mixed volume (with Minkowski addition) and positive products (with addition of divisors).  We show that in some ways this correspondence works well, but in other ways it fails.
The next proposition corresponds to the fact that the function $V(K^{n-1},-)$ is linear on $\mathfrak{K}^{n}$.

\begin{prop}\label{pullback}
Let $X$ be a smooth projective variety of dimension $n$. Assume that $L$ is pseudo-effective and $M$ is movable, then
\begin{equation*}
  \langle L^{n-1}\rangle\cdot M = \langle L^{n-1}\cdot M\rangle.
\end{equation*}

\end{prop}

\begin{proof}
This is an application of \cite[Lemma 6.21]{fl14}. For completeness, we include some arguments here. Firstly, by \cite{nystrom2016duality} the results in \cite{BFJ09} can extend to transcendental classes.

By taking limits, we can assume that $M=\pi_* A$ for some birational modification $\pi$ and K\"ahler class $A$ upstairs. We claim that $\pi^* \langle L^{n-1}\rangle =  \langle (\pi^*L)^{n-1}\rangle$. Note that $\pi_* \langle (\pi^*L)^{n-1}\rangle = \langle L^{n-1}\rangle$, so we need to check that $\langle (\pi^*L)^{n-1}\rangle \cdot E =0$ for any irreducible $\pi$-exceptional divisor $E$. For any $c>0$, the restriction $\pi^*L + cE |_E$ can not be pseudo-effective, thus $\pi^* L$ can not be $E$-big. By \cite[Lemma 4.10]{BFJ09}, we get that $\langle (\pi^*L)^{n-1}\rangle \cdot E =0$, concluding the proof of our claim. Since $\pi^* (\pi_* A) \succeq A$, by basic properties of positive products, we have
\begin{align*}
  \langle L^{n-1} \cdot \pi_* A\rangle &= \langle (\pi^*L)^{n-1} \cdot \pi^*(\pi_* A)\rangle\\
  & \geq \langle (\pi^*L)^{n-1} \cdot A\rangle = \langle (\pi^*L)^{n-1} \rangle\cdot A\\
  & = \pi^* \langle L^{n-1}\rangle \cdot A = \langle L^{n-1}\rangle \cdot \pi_*A
\end{align*}

On the other hand, we always have $\langle L^{n-1}\rangle\cdot M \geq \langle L^{n-1}\cdot M\rangle$. Thus we have the equality $\langle L^{n-1}\rangle\cdot M = \langle L^{n-1}\cdot M\rangle$.
\hfill $\Box$
\end{proof}

\begin{rmk} \rm
In the above proof, the key ingredient is the equality $\pi^*\langle L^{n-1}\rangle = \langle (\pi^*L)^{n-1}\rangle$. An alternative way to see this equality is to use the function $\mathfrak{M}$ studied in \cite{xiao15, lehmann2016positiivty}. For any big class $N$, we have
\begin{align*}
  \frac{\pi^*\langle L^{n-1}\rangle \cdot N}{\vol(N)^{1/n}} &\geq \frac{\langle L^{n-1}\rangle \cdot \pi_* N}{\vol(\pi_*N)^{1/n}}.
\end{align*}
Thus the infimum of LHS is obtained by some big class $N$ satisfying $\pi_* N =L$, which implies that the infimum of LHS is obtained by the positive part of $\pi^*L$. By \cite{lehmann2016positiivty} we get that $$\pi^*\langle L^{n-1}\rangle = \langle P(\pi^*L)^{n-1}\rangle = \langle (\pi^*L)^{n-1}\rangle.$$
\end{rmk}

\begin{cor}
Let $X$ be a smooth projective variety of dimension $n$. Assume that $L$ is pseudo-effective and $M_1, M_2$ are movable, then
\begin{equation*}
  \langle L^{n-1}\cdot (M_1 +M_2)\rangle = \langle L^{n-1}\cdot M_1\rangle + \langle L^{n-1}\cdot M_2\rangle.
\end{equation*}
\end{cor}

However, the more general statement $\langle L_1 \cdot L_2 \cdot ...\cdot L_{n-1}\rangle \cdot M = \langle L_1 \cdot L_2 \cdot ...\cdot L_{n-1}\cdot M\rangle$ is false.  Equivalently, it is not true that $\phi^{*}\langle L_1 \cdot L_2 \cdot ...\cdot L_{n-1}\rangle = \langle \phi^{*} L_1 \cdot \phi^{*} L_2 \cdot ...\cdot \phi^{*} L_{n-1}\rangle$ for a birational map $\phi$.  If this property were true, it would easily imply that the volume function is polynomial (not just locally polynomial) everywhere in the interior of the movable cone.  But this property fails in the following examples suggested by Lazarsfeld and Totaro.

\begin{exmple}[A flop] \rm
Let $X$ be the projective bundle (of quotients) over $\mathbb{P}^{1}$ defined by $\mathcal{O} \oplus \mathcal{O} \oplus \mathcal{O}(-1)$.  There are two natural divisor classes on $X$: the class $f$ of the fibers of the projective bundle and the class $\xi$ of the sheaf $\mathcal{O}_{X/\mathbb{P}^{1}}(1)$.  Using for example \cite{fulger11} and \cite{fl14}, one sees that $f$ and $\xi$ generate the algebraic cohomology classes with the relations $f^{2} = 0$, $\xi^{2}f = -\xi^{3} = 1$ and
\begin{equation*}
\Eff^{1}(X) = \Mov^{1}(X) = \langle f, \xi \rangle \qquad \qquad \Nef^{1}(X) = \langle f, \xi + f \rangle
\end{equation*}
and
\begin{align*}
\Eff_{1}(X) = \langle \xi f, \xi^{2} \rangle \qquad & \qquad \Nef_{1}(X) = \langle \xi f, \xi^{2} + \xi f \rangle.
\end{align*}

We can calculate volume by pushing forward.  We have if $b \geq a$,
\begin{align*}
H^{0}(X, a \xi + b f) & = H^{0}(\mathbb{P}^{1}, \Sym^{a}(\mathcal{O} \oplus \mathcal{O} \oplus \mathcal{O}(-1)) \otimes \mathcal{O}(b)) \\
& = (b+1)(a+1) + ba + (b-1)(a-1) + \ldots + (b-a)(a-a) \\
& = ba^{2}/2 - a^{3}/6 + 3ab/2 + b + 7a/6 +1
\end{align*}
and if $b \leq a$,
\begin{align*}
H^{0}(X, a \xi + b f) &  = (b+1)(a+1) + ba + (b-1)(a-1) + \ldots + (b-b)(a-b) \\
& = ab^{2}/2 - b^{3}/6 + 3ba/2 + a + 7b/6 +1
\end{align*}

So for $b \geq a$ the volume is $3ba^{2} - a^{3}$, and for $a \geq b$ it is $3ab^{2} - b^{3}$.  Note that in a neighborhood of $a=b$ the volume is not locally polynomial.

\end{exmple}

\begin{exmple}[Projective bundle over an abelian surface] \rm
In fact, the volume function can even be transcendental on certain regions of the big movable cone.
\cite{klm13} gives an example of a $\mathbb{P}^{2}$-bundle over an abelian twofold whose volume function is transcendental on a region of the big cone of divisors.  While this region only contains non-movable divisors, their computation also extends into the movable cone as follows.

As in \cite[Section 3]{klm13} let $A$ be the self-product of an elliptic curve without CM and let $f_{1}$, $f_{2}$, and $\Delta$ denote the fibers of the two projection maps and the diagonal.  Consider the projective bundle $\pi: X \to A$ defined as
\begin{equation*}
X = \mathbb{P}_{A}(\mathcal{O}(f_{1}+f_{2}+\Delta) \oplus \mathcal{O}(-f_{1}) \oplus \mathcal{O}(-f_{2}))
\end{equation*}
and let $L$ be a divisor representing $\mathcal{O}_{X/A}(1)$.

For positive real numbers $c_{1},c_{2},c_{3}$ \cite[Section 3]{klm13} computes the volume of $L + \pi^{*}(c_{1}f_{1} + c_{2}f_{2} + c_{3}\Delta)$ as an explicit integral over a region in the plane.  This divisor is movable but not nef when $c_{1} = \frac{1}{2}+\epsilon, c_{2} = \frac{1}{2}+\epsilon, c_{3}=1$ for $\epsilon$ sufficiently small (for example, $\epsilon = 1/10$ will work).  Following the lengthy computation of \cite{klm13}, we find that for rational values of the $c_{i}$ in a neighborhood of the divisor $L + \frac{6}{10}\pi^{*}f_{1} + \frac{6}{10}\pi^{*}f_{2} + \pi^{*}\Delta$ the volume can be written as a (non-zero) quadratic irrational multiple of $\mathrm{ln}$ of a (non-constant) quadratic irrational, plus an algebraic number.

\end{exmple}

\begin{rmk} \rm
In the algebraic setting, for any big and movable classes $L_{1},\ldots,L_{n-1}$ the difference $\phi^{*}\langle L_1 \cdot L_2 \cdot ...\cdot L_{n-1}\rangle - \langle \phi^{*} L_1 \cdot \phi^{*} L_2 \cdot ...\cdot \phi^{*} L_{n-1}\rangle$ is a pseudo-effective $\phi$-contracted curve class.  This follows from the argument establishing negativity of contraction.
\end{rmk}

\section{Alexandrov body and Zariski decompositions}
\label{section alex body}

In this section, we interpret a classical construction due to Alexandrov \cite{alexandrov1938} (an English translation is available in his selected works \cite{alexandrovbook}) as a Zariski decomposition structure in convex geometry.  This is the convex geometry version of the $\sigma$-decomposition (or divisorial Zariski decomposition) of divisors in algebraic geometry.

\subsection{Zariski decomposition}
First we recall the notion of a $\sigma$-decomposition.  Let $X$ be a smooth projective variety of dimension $n$, and let $L\in \Eff^1 (X)$ be a big divisor class.   As in the theories developed in
\cite{Bou04, Nak04, BFJ09, BDPP13}, there exists a decomposition
\begin{equation*}
  L=P(L)+N(L)
\end{equation*}
where $P(L) = \langle L \rangle$ and $N(L) = L - P(L)$.  This decomposition is known as the $\sigma$-decomposition; it should be interpreted as replacing a pseudo-effective class $L$ with a movable class $P(L) \leq L$ with no ``loss'' of positivity.  The key properties of this decomposition are the following:
\begin{enumerate}
\item $P(L)\in \Mov^1 (X)$ and $N(L)$ is a rigidly embedded effective divisor class\footnote{Here we use the terminology of \cite{Bou04}, ``rigidly embedded'' means that the class $N(L)$ contains only one positive current.}.
\item $\vol(L)=\vol(P(L))$.
\item The orthogonality property: $\langle P(L) ^{n-1}\rangle \cdot N(L) =0$.
\end{enumerate}
The first two properties characterize the $\sigma$-decomposition: if $L = P+N$ is any decomposition such that $P\in \Mov^1 (X), N\in \Eff^1 (X)$ and $\vol(L)=\vol(P)$, then it is shown in \cite[Section 5]{fl14} that we must have $P= P(L), N= N(L)$.

In fact, \cite{lehmann2016positiivty} shows a stronger property: the $\sigma$-decomposition exactly captures the failure of the volume function to be strictly log concave.

\begin{thrm}[see \cite{lehmann2016positiivty} Theorem 1.6]
Let  $X$ be a smooth projective variety of dimension $n$.  Let $L$ and $D$ be big divisor classes on $X$.  Then we have
$$\vol(L + D)^{1/n} = \vol(L)^{1/n} + \vol(D)^{1/n}$$
if and only if $P(L), P(D)$ are proportional.
\end{thrm}

\subsection{Alexandrov body}
Let us recall the construction of the Alexandrov body (see e.g. \cite[Section 7.5]{schneiderconvex}). Let $\mathfrak{K}_{0}^n \subseteq \mathfrak{K}^n$ be the space of convex bodies with $0$ as an interior point. Let $f$ be a strictly positive continuous function.  The Alexandrov body associated with $f$ is defined as the unique maximal element (denoted by $K$), with respect to set inclusion, of
the set
\begin{equation*}
  \{Q\in \mathfrak{K}_{0}^n|\ h_Q \leq f \}.
\end{equation*}
Equivalently, this means that the support function $h_K = \sup \rho$ where the sup is taken over all 1-homogeneous convex functions less than or equal to $f$.  The convex body $K$ can also be characterized as follows:
\begin{equation*}
  K=\bigcap_{y\in \mathbb{S}^{n-1}} \{x\in \mathbb{R}^n |\ x\cdot y \leq f(y)\}.
\end{equation*}
If we denote $P(f):= h_K$, then we have
\begin{equation}\label{eq othogonality alex1}
  P(f) = f,\ a.e.\ \textrm{with respect to}\ \mu_K,
\end{equation}
where $\mu_K :=S(K^{n-1};\cdot)$ is the area measure of $K$.
Equivalently, we get
\begin{equation}\label{eq othogonality alex2}
  (f-P(f))\cdot \mu_K :=\int_{\mathbb{S}^{n-1}}(f-P(f))d \mu_K=0.
\end{equation}
The volume of $f$ is defined to be
\begin{equation}\label{eq vol alex}
  \vol(f):=\vol(K).
\end{equation}
For $L\in \mathfrak{K}_0 ^n$, it is clear that $P(h_L)=h_L, \vol(h_L)=\vol(L)$.

\subsubsection{Decomposition theory}
As an immediate consequence of the characterization of Alexandrov body, we have

\begin{lem}
Let $f, g$ be two strictly positive continuous functions such that $f=g +l_a$ for some linear function $l_a$, then we have $P(f)=P(g)+l_a$.
\end{lem}

This implies the decomposition of $\mathcal{C}_p$:

\begin{thrm}[Decomposition of $(\mathcal{C}_p, \vol)$ \cite{alexandrov1938}] \label{thrm alexzar}
Assume that $f\in \mathcal{C}_p$ is represented by a strictly positive continuous function.  Then there exists a decomposition
\begin{equation*}
  f=P(f)+N(f)
\end{equation*}
satisfying:
\begin{itemize}
  \item $P(f)\in \mathcal{C}$ and $N(f)=f-P(f)\in \mathcal{C}_p$;
  \item $\vol(f)=\vol(P(f))$;
  \item if we assume $P(f)$ is represented by a convex body $K$ and denote $\langle P(f)^{n-1}\rangle:= S(K^{n-1};\cdot)$, then the orthogonality estimate holds in the sense that $\langle P(f)^{n-1}\rangle \cdot N(f) =0$.
\end{itemize}
\end{thrm}

We can rephrase this result to look exactly analogous to the characterization of the $\sigma$-decomposition described above; by the Brunn-Minkowski theorem we obtain

\begin{prop}
Suppose that $f,g \in \mathcal{C}_{p}$ are represented by strictly positive continuous functions.  Then we have $$\vol(f+g)^{1/n} = \vol(f)^{1/n} + \vol(g)^{1/n}$$
if and only if $P(f), P(g)$ are proportional.
In particular, if $f = P+N$ is a decomposition such that $P\in \mathcal{C}, N\in \mathcal{C}_p$ and $\vol(P)=\vol(f)$, then we must have $P=P(f), N= N(f)$.
\end{prop}

\begin{proof}
The proof is similar to the complex geometry setting. We only need to prove the second point. Integrating $f= P+N$ against with the area measure $\langle P(f)^{n-1}\rangle$ implies
\begin{align*}
  \vol(f) &= P \cdot \langle P(f)^{n-1}\rangle + N\cdot \langle P(f)^{n-1}\rangle\\
  & \geq P \cdot \langle P(f)^{n-1}\rangle \geq \vol(P)^{1/n} \vol(P(f))^{n-1/n}\\
  & = \vol(f),
\end{align*}
where the first inequality follows from $N \in \mathcal{C}_p$, the second inequality follows from the Brunn-Minkowski inequality and the last equality follows from the assumption $\vol(P)=\vol(f)$. Thus we get equality everywhere, in particular, $P \cdot \langle P(f)^{n-1}\rangle = \vol(P)^{1/n} \vol(P(f))^{n-1/n}$. Note that $P \cdot \langle P(f)^{n-1}\rangle$ is just the mixed volume computed by the convex bodies representing $P$ and $P(f)$. By Brunn-Minkowski theory, up to some translation these two convex bodies are equal, which in turn implies $P = P(f)$.
\hfill $\Box$
\end{proof}

\begin{rmk} \rm
Theorem \ref{thrm alexzar} can also be seen as some kind of infinite dimensional extension of the theory developed in \cite{lehmxiao16convexity} when applied to this particular setting.
\end{rmk}

\begin{rmk} \rm
One should compare the Alexandrov body construction with the Berman-Boucksom orthogonality estimate in complex pluripotential theory (see \cite{berman2010growth}), where the authors apply a complex envelope construction analogous to the construction of $P(f)$. This complex envelope construction admits similar orthogonality properties as the Zariski decomposition, and has played an important role in the variational approach to solve complex Monge-Amp\`ere equations (see \cite{BBGZvariational}). Indeed, by the theory developed in \cite{boucksom2015solution} on non-Archimedean Monge-Amp\`ere equations, the orthogonality in the divisorial Zariski decomposition is closely related to non-Archimedean Berman-Boucksom orthogonality.
\end{rmk}

\section{Morse type inequalities} \label{morseineq section}
Let $X$ be a compact K\"ahler manifold of dimension $n$, and let $A, B\in H^{1,1}(X, \mathbb{R})$ be nef classes.   Demailly's conjecture on the (weak) transcendental Morse inequality says that $$\vol(A-B)\geq A^n -n A^{n-1}\cdot B.$$  This important conjecture was recently proved for projective manifolds in \cite{nystrom2016duality}.  For non-projective K\"ahler manifolds, it is proved in \cite{popovici2016sufficientbig} that $A^n -n A^{n-1}\cdot B>0$ implies $A-B$ is big (see also \cite{xiao13weakmorse} for a weaker result).

As an easy corollary of the main result of \cite{nystrom2016duality}, we have

\begin{prop}
Let $X$ be a smooth projective variety of dimension $n$. If $L, M\in H^{1,1} (X, \mathbb{R})$ are movable classes, then
\begin{equation*}
  \vol(L-M)\geq \vol(L)-n\langle L^{n-1}\rangle\cdot M.
\end{equation*}
\end{prop}

\begin{proof}
For reader's convenience, we repeat some arguments from \cite{xiao2014movable}. Without loss of generality, we can assume that $L$ is big and that $M=\mu_* \widehat{\omega}$ for a birational modification $\mu$ and for some K\"ahler class $\widehat{\omega}$ upstairs. Take a suitable sequence of Fujita approximations, $\pi_m ^*(\mu^* L) = \omega_m +[E_m]$, such that
\begin{equation*}
  \vol(\mu^* L) = \lim_{m\rightarrow \infty} \omega_m ^n, \ \ \ \langle (\mu^* L)^{n-1}\rangle\cdot \widehat{\omega} = \lim_{m\rightarrow \infty} \omega_m ^{n-1} \cdot \pi_m ^* \widehat{\omega}.
\end{equation*}

By \cite{nystrom2016duality}, we always have $\vol(\omega_m - \pi_m ^* \widehat{\omega}) \geq \omega_m ^n -n \omega_m ^{n-1} \cdot \pi_m ^* \widehat{\omega}$, which implies that
\begin{align*}
  \vol(\mu^* L - \widehat{\omega}) &= \vol(\pi_m ^* (\mu^* L - \widehat{\omega})) \\
  &\geq \vol(\omega_m - \pi_m ^* \widehat{\omega})\\
  &\geq \omega_m ^n -n \omega_m ^{n-1} \cdot \pi_m ^* \widehat{\omega}.
\end{align*}
Let $m$ tend to infinity, we get $\vol(\mu^* L - \widehat{\omega}) \geq \vol(L) - n \langle (\mu^* L)^{n-1}\rangle\cdot \widehat{\omega}$. This implies the desired inequality
\begin{equation*}
  \vol(L - \mu_*\widehat{\omega}) \geq \vol(L) - n \langle L^{n-1}\rangle\cdot \mu_*\widehat{\omega}.
\end{equation*}
\hfill $\Box$
\end{proof}

\subsection{Convex geometry analogue}\label{morse convexgeom}
Due to the importance of the Morse inequality in complex geometry, it is interesting to ask for an analogue for convex bodies.  More precisely, let $K, L \in \mathfrak{K}^n$ and let $h_K, h_L\in \mathcal{C}$ be their support functions, one may wonder if $$\vol(h_K -h_L)\geq \vol(K)-nV(K^{n-1}, L)$$ holds, and in particular, $h_K -h_L$ could be represented by a strictly positive function if $\vol(K)-nV(K^{n-1}, L)>0$.

Note that we can assume $\vol(K)-nV(K^{n-1}, L) > 0$, since the inequality is trivial otherwise.  It is also easy to see that the inequality holds if $K, L$ are proportional up to some translation, so we may assume otherwise.

We first verify that when $\vol(K)-nV(K^{n-1}, L)>0$, the function $h_{K}-h_{L}$ can be represented by a strictly positive continuous function.  This can be accomplished by applying the Diskant inequality (see \cite[Section 7.2]{schneiderconvex}). Recall that the \emph{inradius} of $K$ relative to $L$ is defined by:
$$r(K, L):=\max \{\lambda>0|\lambda L +t \subset K\ \textrm{for some}\ t\in \mathbb{R}^n\}. $$
Applying the Diskant inequality to $K, L$, we get
\begin{align*}
  r(K, L)&\geq \frac{V(K^{n-1}, L)^{1/n-1}-\left(V(K^{n-1}, L)^{n/n-1}-\vol(K)\vol(L)^{1/n-1}\right)^{1/n}}{\vol(L)^{1/n-1}}.
  \end{align*}
The generalized binomial formula yields
\begin{align*}
  &V(K^{n-1}, L)^{1/n-1}-\left(V(K^{n-1}, L)^{n/n-1}-\vol(K)\vol(L)^{1/n-1}\right)^{1/n}\\
  &= \sum_{k=1}^{\infty} \left( \begin{array}{c} 1/n \\ k \end{array}\right) (-1)^{k+1} \vol(K)^k\vol(L)^{k/n-1} V(K^{n-1}, L)^{\frac{1-kn}{n-1}}\\
  &= \frac{1}{n}\vol(K)\vol(L)^{1/n-1} V(K^{n-1}, L)^{-1} + \sum_{k=2}^{\infty}...\\
  & \geq \frac{1}{n}\vol(K)\vol(L)^{1/n-1} V(K^{n-1}, L)^{-1},
\end{align*}
because every term in the sum $\sum_{k=2}^{\infty} ...$ is non negative. Substituting, we obtain
\begin{equation*}
 r(K,L) \geq \frac{\vol(K)}{nV(K^{n-1}, L)}.
\end{equation*}

Now the assumption $\vol(K)-nV(K^{n-1}, L)>0$ implies that, up to a translation, $K\supset (1+\delta)L$ for some $\delta>0$. In particular, up to a linear function, we have $h_K -h_L \geq \delta h_L$. This implies that the class $\{h_K -h_L \}\in \mathcal{C}_p$ can be represented by a strictly positive continuous function.

\begin{rmk} \rm
It is of course not true in general that $h_{K} - h_{L}$ can be represented by a strictly positive convex continuous function.  For example, this necessarily fails when $K$ is indecomposable as discussed in Section \hyperref[indecomposability sec]{5.A.a}.
\end{rmk}

Now the Morse inequality for convex bodies follows from the following result due to \cite{alexandrov1938} (see also \cite[Section 7.5]{schneiderconvex}). 

\begin{lem}
Let $f\in \mathcal{C}_p$ be a strictly positive continuous function with its Alexandrov body given by $P(f)=h_Q$.  For any continuous function $g\in V$ we have
\begin{equation*}
   \left. \frac{d}{dt} \right|_{t = 0} \vol(f+tg)=n \langle P(f)^{n-1}\rangle \cdot g.
\end{equation*}
\end{lem}

\noindent Applying this result to our original set-up, we see that for $t\in [0, 1]$
\begin{align*}
  \vol(h_K -h_L)-\vol(K) &=\int_0 ^1 \frac{d}{dt} \vol(h_K -t h_L)dt\\
  &=\int_0 ^1 -n \langle P(h_K -t h_L)^{n-1}\rangle \cdot h_L dt\\
  &\geq -n V(K^{n-1}, L),
\end{align*}
where the last inequality follows from the monotone property of mixed volumes.

Summarizing the discussions above, we get
\begin{thrm}[Morse inequality for convex bodies]
\label{thrm morseconvx}
Let $h_K, h_L \in \mathcal{C}$.  If
$$\vol(K)-nV(K^{n-1}, L)>0,$$
then $h_K -h_L \in \mathcal{C}_p ^{\circ}$ which means that it can be represented by some strictly positive continuous function. Moreover, we have $\vol(h_K -h_L)\geq \vol(K)-nV(K^{n-1}, L)$.
\end{thrm}

For convex bodies, we also have the ``reverse'' Khovanskii-Teissier inequality (see \cite[Remark 4.18]{lehmxiao16convexity} for this terminology).

\begin{prop}
\label{reverse KT convex1}
Fix three convex bodies $K, L, M\in \mathfrak{K}^n$, then
\begin{align*}
  nV(L, K^{n-1})V(K, M^{n-1})\geq \vol(K)V(L, M^{n-1}).
\end{align*}
\end{prop}

\begin{proof}
By taking limits, we can assume the convex bodies have non-empty interior. After scaling, we can assume $V(K, M^{n-1})=V(L, M^{n-1})$. Under this assumption, we need to show $nV(L, K^{n-1})\geq \vol(K)$.  Otherwise we would have $h_{K} - h_{L} \in \mathcal{C}_{p}^{\circ}$ by Theorem \ref{thrm morseconvx}.  Thus $$\int_{\mathbb{S}^{n-1}}(h_K-h_L)dS(M^{n-1})=V(K, M^{n-1})-V(L, M^{n-1})>0,$$
a contradiction. \hfill $\Box$
\end{proof}

\subsubsection{Indecomposability} \label{indecomposability sec}

It is interesting to ask when a big movable $(1,1)$-class lies on the boundary (or extremal ray) of the movable cone.  The analogue in the convex geometry setting is the notion of indecomposability.  Recall that $K$ is \emph{decomposable} if there exist convex bodies $M, D$ not homothetic to $K$ such that $K=M+D$. Otherwise, we call that $K$ is \emph{indecomposable}. In particular, the existence of indecomposable convex bodies implies that, in the above Morse type inequality for convex bodies, in general $h_K -h_L$ cannot be represented by a support function.
We refer the readers to \cite[Chapter 3]{schneiderconvex} for decomposable and indecomposable convex bodies.

\begin{exmple} \rm
Suppose that $X$ is a toric variety defined by a fan $\Sigma$.  Suppose that $P$ is an indecomposable polytope whose normal fan is refined by $\Sigma$.  Then $P$ induces a big movable divisor on the boundary of the movable cone.  The converse is false -- however, if one focuses only on integer multiples of lattice polytopes, there is a good converse statement as in \cite[Proposition 6.2.13]{cls2011toric}.

For example, the only two-dimensional indecomposable convex bodies are triangles.  In the toric setting, a triangle represents the pullback under a birational contraction of an ample generator of a ($\mathbb{Q}$-factorial) toric surface of Picard rank $1$.  Such a divisor is necessarily on the boundary of the movable cone.
\end{exmple}

\subsection{Reverse Khovanskii-Teissier inequality}
In the algebraic setting, an important step in the analysis of the Morse type bigness criterion is the ``reverse" Khovanskii-Teissier inequality for big and nef divisors $A$, $B$, and a movable curve class $\beta$:
\begin{equation*}
n(A \cdot B^{n-1})(B \cdot \beta) \geq B^{n} (A \cdot \beta).
\end{equation*}
We prove a more general statement on ``reverse" Khovanskii-Teissier inequalities in the analytic setting\footnote{Indeed, this kind of generalization has already been noted in \cite[Remark 3.1]{xiao13weakmorse}, and it was also noted that the constant ${k!(n-k)!}/{4n!}$ there could be improved by combining with the technique of \cite{popovici2016sufficientbig}.
Some related work has also appeared independently in the recent preprint \cite{popovici15}.}.

\begin{thrm}
\label{thrm appendix reserve KT}
Let $X$ be a compact K\"ahler manifold of dimension $n$. Let $\omega, \beta, \gamma$ be three nef $(1,1)$ classes on $X$. Then we have
$$
(\beta^k \cdot \alpha^{n-k})\cdot (\alpha^k\cdot \gamma^{n-k})\geq \frac{k!(n-k)!}{n!}\alpha^n\cdot (\beta^k \cdot\gamma^{n-k}).
$$
\end{thrm}

\begin{proof}
The proof depends on solving Monge-Amp\`{e}re equations and the method of \cite{popovici2016sufficientbig}. Without loss of generality, we can assume that $\gamma$ is normalised such that $\beta^k\cdot \gamma^{n-k}=1$. Then we need to show
\begin{align}
\label{eq need to prove}
(\beta^k \cdot \alpha^{n-k})\cdot (\alpha^k\cdot \gamma^{n-k})\geq \frac{k!(n-k)!}{n!}\alpha^n.
\end{align}
We first assume $\alpha,\beta,\gamma$ are all K\"ahler classes. We will use the same symbols to denote the K\"ahler metrics in corresponding K\"ahler classes. By the Calabi-Yau theorem \cite{Yau78}, we can solve the following Monge-Amp\`{e}re equation:
\begin{align}
\label{eq MA}
(\alpha+i\partial\bar \partial \psi)^n= \left(\int \alpha^n \right) \beta^k\wedge \gamma^{n-k}.
\end{align}
Denote by $\alpha_\psi$ the K\"ahler metric $\alpha+i\partial\bar \partial \psi$. Then we have
\begin{align*}
(\beta^k \cdot \alpha^{n-k})\cdot (\alpha^k\cdot \gamma^{n-k})&=\int \beta^k \wedge \alpha_\psi^{n-k}\cdot \int \alpha_\psi^k \wedge \gamma^{n-k}\\
&= \int \frac{\beta^k \wedge \alpha_\psi^{n-k}}{\alpha_\psi ^n}\alpha_\psi ^n\cdot \int \frac{\alpha_\psi^k \wedge \gamma^{n-k}}{\alpha_\psi ^n}\alpha_\psi ^n\\
&\geq \left( \int \left(\frac{\beta^k \wedge \alpha_\psi^{n-k}}{\alpha_\psi ^n} \cdot \frac{\alpha_\psi^k \wedge \gamma^{n-k}}{\alpha_\psi ^n} \right)^{1/2}\alpha_\psi ^n\right)^2.
\end{align*}

The last line follows because of the Cauchy-Schwarz inequality. We claim that the following pointwise inequality holds:
\begin{align}
\label{eq pointwise}
\frac{\beta^k \wedge \alpha_\psi^{n-k}}{\alpha_\psi ^n} \cdot \frac{\alpha_\psi^k \wedge \gamma^{n-k}}{\beta^k \wedge \gamma^{n-k}} \geq \frac{k!(n-k)!}{n!}.
\end{align}
Then by (\ref{eq MA}) it is clear that the above pointwise inequality implies the desired inequality (\ref{eq need to prove}).
Now we prove the claim. For any fixed point $p\in X$, we can choose some coordinates such that at the point $p$: $$\alpha_\psi= i\sum_{j=1}^n dz^j \wedge d\bar z ^j,
\quad \beta = i\sum_{j=1}^n \mu_j dz^j \wedge d\bar z ^j,
$$
and
 $$\gamma^{n-k}=i^{n-k} \sum_{|I|=|J|=n-k} \Gamma_{IJ} dz_I \wedge d\bar z_{J}.$$
Denote by $\mu_J$ the product $\mu_{j_1}...\mu_{j_k}$ with index $J=(j_1<...<j_k)$ and denote by $J^c$ the complement index of $J$. Then it is easy to see at the point $p$ we have
$$
\frac{\beta^k \wedge \alpha_\psi^{n-k}}{\alpha_\psi ^n} \cdot \frac{\alpha_\psi^k \wedge \gamma^{n-k}}{\beta^k \wedge \gamma^{n-k}}= \frac{k!(n-k)!}{n!} \frac{(\sum_{J}\mu_J)(\sum_K \Gamma_{KK})}{\sum_{J}\mu_J \Gamma_{J^c J^c}}\geq \frac{k!(n-k)!}{n!}.
$$
This then finishes the proof of the case when $\alpha,\beta,\gamma$ are all K\"ahler classes. If they are just nef classes, by taking limits, then we get the desired inequality.
\hfill $\Box$
\end{proof}

By taking suitable Fujita approximations, we immediately get

\begin{cor}
\label{psef reserve KT}
Let $X$ be a compact K\"ahler manifold of dimension $n$. Let $\alpha, \beta, \gamma$ be three pseudo-effective $(1,1)$ classes on $X$. Then we have
$$
\langle\beta^k \cdot \alpha^{n-k}\rangle\cdot \langle\alpha^k\cdot \gamma^{n-k}\rangle\geq \frac{k!(n-k)!}{n!}\vol(\alpha) \langle\beta^k \cdot\gamma^{n-k}\rangle.
$$
\end{cor}

\begin{proof}
By the continuity of positive products, we can assume that the classes are big. We may take a sequence of Fujita approximations $\mu_m : X_m \rightarrow X$ such that
\begin{equation*}
  \mu_m ^* \alpha = \alpha_m + [E_m],\ \ \mu_m ^* \beta = \beta_m + [F_m],\ \ \mu_m ^* \gamma = \gamma_m + [G_m],
\end{equation*}
where $\alpha_m, \beta_m, \gamma_m$ are K\"ahler and $E_m, F_m, G_m$ are effective. Moreover, the positive products are the limits of intersections of K\"ahler classes upstairs, for instance, $$\langle \beta^k \cdot \alpha^{n-k}\rangle= \lim_{m\rightarrow \infty} \beta_m ^k \cdot \alpha_m ^{n-k}.$$
Then the inequality follows from Theorem \ref{thrm appendix reserve KT} by taking limits.
\hfill $\Box$
\end{proof}

In the convex geometry setting, it is natural to ask whether an analogue holds for convex bodies. As the generalization of Proposition \ref{reverse KT convex1},
we get the following result.

\begin{thrm}
Assume that $K, L, M$ are three convex bodies in $\mathbb{R}^n$.  Then
\begin{align*}
V(K^k, L^{n-k})V(L^k, M^{n-k})\geq \frac{k!(n-k)!}{n!}\vol(L)V(K^k, M^{n-k}).
\end{align*}
\end{thrm}

Since we did not find a similar statement in the literature, we provide a brief proof below. The proof is inspired by \cite{milman99masstransport}, where the authors reproved some of the Alexandov-Fenchel inequalities by using mass transport (or Brenier maps). In our setting, instead of using complex Monge-Amp\`{e}re equations for $(1,1)$ classes as in Theorem \ref{thrm appendix reserve KT}, we apply mass transport (real Monge-Amp\`{e}re equations) to convex bodies. Then we reduce the inequality to an inequality for mixed discriminants, which is easy to prove.

\begin{proof}
By the continuity of mixed volumes and taking a limit, it suffices to consider the case when $K,L,M$ have non-empty interior.  We can also assume that the convex bodies are open. Without loss of generality, we can assume $\vol(L)=1$.

By \cite{gromov1990convex}, for $K, M$, we can take two $C^2$ strictly convex functions $f_K, f_M: \mathbb{R}^n \rightarrow \mathbb{R}$ such that
\begin{equation*}
  \Image(\nabla f_K)= K,\ \ \  \Image(\nabla f_M)=M
\end{equation*}
where $\nabla$ is the gradient operator. In the following, we denote the Hessian operator by $\nabla^2$.

For $n\times n$ positive matrices $M_1, ..., M_r$, the discriminant $D(M_1 ^{k_1},...,M_r ^{k_r})$ is given by the following expansion:
\begin{equation*}
  \det(\sum_{k=1} ^r t_k M_k)=\sum_{k_1 +...+k_r=n}  \frac{n!}{k_1 !...k_r !}D(M_1 ^{k_1},...,M_r ^{k_r}) t_1 ^{k_1}...t_r ^{k_r}.
\end{equation*}

Let
\begin{equation*}
  c=\int_{\mathbb{R}^n} D((\nabla^2 f_K)^k, (\nabla^2 f_M)^{n-k})dx.
\end{equation*}
The measure $\rho dx=c^{-1}D((\nabla^2 f_K)^k, (\nabla^2 f_M)^{n-k})dx$ is a probability measure on $\mathbb{R}^n$. Note that by the argument in \cite{milman99masstransport} we have
\begin{equation*}
  V(K^k, M^{n-k})=c.
\end{equation*}

By \cite{brenier1991polar} or \cite{mccann95existence}, we can take a convex function $F_L$ satisfying $\Image(\nabla F_L) =L$ and $\det(\nabla^2 F_L) = \rho$. As the argument in \cite{milman99masstransport} again, we have
\begin{align*}
  V(K^k, L^{n-k})&=\int_{\mathbb{R}^n} D((\nabla^2 f_K)^k, (\nabla^2 F_L)^{n-k})dx, \\
  V(L^k, M^{n-k})&=\int_{\mathbb{R}^n} D((\nabla^2 F_L)^k, (\nabla^2 f_M)^{n-k})dx.
\end{align*}

As in the proof of Theorem \ref{thrm appendix reserve KT}, by Cauchy-Schwarz inequality we get
\begin{align*}
   V(K^k, L^{n-k})V(L^k, M^{n-k})&=\int_{\mathbb{R}^n} D((\nabla^2 f_K)^k, (\nabla^2 F_L)^{n-k})dx \int_{\mathbb{R}^n} D((\nabla^2 F_L)^k, (\nabla^2 f_M)^{n-k})dx\\
   &\geq \left(\int_{\mathbb{R}^n} (D((\nabla^2 f_K)^k, (\nabla^2 F_L)^{n-k})D((\nabla^2 F_L)^k, (\nabla^2 f_M)^{n-k}))^{1/2}dx\right)^2.
\end{align*}
Similar to the pointwise inequality (\ref{eq pointwise}), for the mixed discriminants we have
\begin{align*}
  D((\nabla^2 f_K)^k, (\nabla^2 F_L)^{n-k}) D((\nabla^2 F_L)^k, (\nabla^2 f_M)^{n-k}) \geq
  &\frac{k!(n-k)!}{n!} \det(\nabla^2 F_L) D((\nabla^2 f_K)^k, (\nabla^2 f_M)^{n-k}).
\end{align*}
This implies that
\begin{align*}
   V(K^k, L^{n-k})V(L^k, M^{n-k})&\geq \left(\int_{\mathbb{R}^n} (\frac{k!(n-k)!}{n!} c^{-1} D((\nabla^2 f_K)^k, (\nabla^2 f_M)^{n-k})^2)^{1/2}dx\right)^2\\
   &=\frac{k!(n-k)!}{n!}V(K^k, M^{n-k})\\
   &=\frac{k!(n-k)!}{n!}\vol(L) V(K^k, M^{n-k}),
\end{align*}
where the last equality follows because we assume $\vol(L)=1$. This finishes the proof of the desired inequality.
\hfill $\Box$
\end{proof}

\section{Cone dualities and Minkowski's existence theorem}\label{section minkowski}
In this section, we discuss the convex geometry analogue of the duality of positive cones in \cite{BDPP13}.

\subsection{The dual cone $\mathcal{C}_p ^*$}

For a projective manifold $X$,  the movable cone $\Mov_1 (X)$ has the following characterization:

\begin{thrm}[see \cite{lehmann2016positiivty}] \label{thrm LX movcone}
Let $X$ be a smooth projective variety of dimension $n$. Then $\Mov_1 (X)$ is the closure of the set $S=\{\langle L^{n-1}\rangle|\ L\ \textrm{big and movable class}\}$. Moreover, the set $S$ is convex.
\end{thrm}

In summary, by \cite{BDPP13} and \cite{lehmann2016positiivty}:
\begin{align*}
  \Mov_1 (X) &=\textrm{the closed cone \emph{genererated} by}\ \{\langle L_1 \cdot ...\cdot L_{n-1}\rangle|\ L_i\ \textrm{big (1,1)-class}\rangle\}\
  \textrm{by \cite{BDPP13, nystrom2016duality}}\\
  & =\textrm{the closure of \emph{the set}}\ \{\langle L^{n-1}\rangle|\ L\ \textrm{big and movable class}\}\ \textrm{by \cite{lehmann2016positiivty}}.
\end{align*}

Motivated by the correspondence $\mathcal{C}_p \leftrightsquigarrow \Eff^1 (X)$, we will consider the dual $\mathcal{C}_{p}^{*}$ by analogue with the duality of cones $\Eff^1 (X)^* = \Mov_1 (X)$.  

The dual space of $C^0 (\mathbb{S}^{n-1})$ is the space of Borel measures on $\mathbb{S}^{n-1}$. Thus a continuous linear functional on $V$ is a Borel measure that vanishes on $L(\mathbb{S}^{n-1})$. Indeed, we have

\begin{lem}\label{lem measure vanlin}
A continuous linear functional $\nu$ belongs to $V^*$ if and only if $\nu$ is a Borel measure that has the origin as its center of mass, that is,
\begin{equation*}
  \int_{\mathbb{S}^{n-1}} x d\nu(x) =0.
\end{equation*}
\end{lem}

The Borel measures in the dual of $\mathcal{C}_{p}$ are precisely the set of nonnegative measures with barycenter $0$.  Alexandrov's version of Minkowski's existence theorem yields an alternative description closer to the complex geometric version above.

To get a feeling for this theorem, let us first recall the original result for polytopes due to Minkowski: let $u_1, ..., u_N$ be pairwise distinct unit vectors linearly spanning $\mathbb{R}^n$, and let $f_1, ..., f_N$ be positive real numbers, then there exists a polytope $P$ (unique up to translation) having precisely $u_1, ..., u_N$ as its outer normal vectors and having $f_1, ..., f_N$ as its corresponded face areas if and only if
\begin{equation*}
  \sum_{i=1} ^N f_i u_i =0.
\end{equation*}
In particular, if we consider $u_i$ as the Dirac measure $\delta_{u_i}$ on the sphere, then $\sum_{i=1} ^N f_i \delta_{u_i}$ is a Borel measure on the sphere with the origin as its center of mass and does not concentrate on any great subsphere. In the general case, we have the following result due to Alexandrov \cite{alexandrov1938} (see also \cite[Chapter V]{alexandrovbook}, \cite[Chapter 8]{schneiderconvex}).

\begin{thrm}[Minkowski's existence theorem]
Assume $\mu$ is a positive Borel measure on $\mathbb{S}^{n-1}$, which is not concentrated on any great subsphere and satisfies
\begin{align*}
\int_{\mathbb{S}^{n-1}} x d\mu(x)=0.
\end{align*}
Then there exists a convex body $K\in \mathfrak{K}^n$, unique up to translation, such that $\mu=S(K,...,K;\cdot)$.
\end{thrm}

Now we can give the description of $\mathcal{C}_{p}^*$, which corresponds exactly to the structure of the movable cone $\Mov_1 (X)$, which follows easily from Minkowski's existence theorem.

\begin{thrm}\label{thrm dualityconvx}
The dual cone $\mathcal{C}_{p}^* \subseteq V^*$ can be characterized as follows:
\begin{align*}
\mathcal{C}_{p}^* &= \textrm{the closed cone \emph{generated} by}\ \{S(K_1, ..., K_{n-1};\cdot)|K_i\in \mathfrak{K}^n\}\\
&=\textrm{the closure of \emph{the set}}\ \{S(K^{n-1};\cdot)|K\in \mathfrak{K}^n\}.
\end{align*}
Alternatively, by the embedding $\Phi$, if we consider the mixed volume $V(K_1, ..., K_{n-1}, \cdot)$ as a continuous linear functional on $V$, then
\begin{align*}
\mathcal{C}_{p}^* &= \textrm{the closed cone \emph{generated} by}\ \{V(K_1, ..., K_{n-1};\cdot)|K_i\in \mathfrak{K}^n\}\\
&=\textrm{the closure of \emph{the set}}\ \{V(K^{n-1};\cdot)|K\in \mathfrak{K}^n\}.
\end{align*}
If we use the notation $\langle h_K ^{n-1}\rangle:= S(K^{n-1};\cdot)$, then we have
\begin{align*}
  \mathcal{C}_{p}^* =\textrm{the closure of \emph{the set}}\ \{\langle P(f)^{n-1}\rangle|f\in \mathcal{C}_p\}.
\end{align*}
\end{thrm}

\begin{proof}
It is obvious that the Borel measures on the right hand side are positive and are in the dual cone $\mathcal{C}_p ^*$. We only need to verify that any element in $\mathcal{C}_{p}^*$ can be written as a limit of a sequence of area measures $S(K, ..., K;\cdot)$ with $K\in \mathfrak{K}^n$.

We claim: let $\nu\in \mathcal{C}_p ^*$, then for any $\epsilon>0$ and $L\in \mathfrak{K}^n$ there exists a convex body $L_\epsilon$ such that
\begin{align*}
\nu+\epsilon S(L^{n-1};\cdot)=S(L_{\epsilon}^{n-1};\cdot).
\end{align*}
To this end, by Lemma \ref{lem measure vanlin} we know that $\nu$ has the origin as its centroid:
\begin{align*}
\int_{\mathbb{S}^{n-1}}xd\nu(x)=0.
\end{align*}

As $\nu$ is a positive measure, for any $L\in \mathfrak{K}^n$ and any $\epsilon>0$, $\nu+\epsilon S(L^{n-1};\cdot)$ does not concentrate on any great subsphere.
Now the existence of $L_\epsilon$ follows immediately from Minkowski's existence theorem.
\hfill $\Box$
\end{proof}

By Minkowski's existence theorem, we have:

\begin{prop}
Assume $K_1,...,K_{n-1}\in \mathfrak{K}^n$ are fixed convex bodies, then there exists a convex body $K \in \mathfrak{K}^n$, unique up to translation, such that $V(K_1,...,K_{n-1},\cdot)=V(K^{n-1},\cdot)$.
\end{prop}

\begin{proof}
We only need to consider the uniqueness, which is an direct corollary of Brunn-Minkowski theory.
\hfill $\Box$
\end{proof}

\begin{rmk} \rm
Let $X$ be a smooth projective variety of dimension $n$, and let $\Eff^1 (X)$ be its pseudo-effective cone of divisor classes. Thanks to the Newton-Okounkov body construction, if we fix a flag $$\mathcal{F}=(X=H_0 \supset H_1 \supset... \supset H_{n-1}\supset H_{n} =\{p\})$$
where $H_i$ is an irreducible smooth subvariety of $\codim =i$, then we have a map $\Delta: \Eff^1 (X)^{\circ} \rightarrow \mathfrak{K}^n$. The map $\Delta$ associates to every big divisor class $L$ a convex body $\Delta(L)$ with non-empty interior satisfying $\vol(L) = \vol(\Delta(L))$. In particular, by the embedding theorem of convex sets, we get a map
\begin{equation*}
  F : \Eff^1 (X)^{\circ} \rightarrow \mathcal{C}.
\end{equation*}
Thus the big cone $\Eff^1 (X)^{\circ}$ can be realized as a finite dimensional convex subcone of $\mathcal{C}$. But $F$ is only concave. In general, for two big divisor classes $L_1, L_2$, we have $\Delta(L_1 +L_2) \supseteq \Delta(L_1) + \Delta(L_2)$, which implies
\begin{equation*}
  F (L_1 +L_2) \succeq F (L_1) +F (L_2).
\end{equation*}
By relating mixed volumes and positive products of big divisors, it would be interesting to see whether Minkowski's existence theorem could provide a new convex geometry proof of the cone duality $\Eff^1 (X)^* =\Mov_1 (X)$.
\end{rmk}

\subsection{Legendre-Fenchel type transforms}

Now that we have the characterization of $\mathcal{C}_p ^*$, we can define a volume function for all strictly positive continuous functions using a Legendre-Fenchel type transform as in \cite{lehmxiao16convexity} (see also \cite{xiao15} for similar volume characterization of divisors.).
Assume that $f\in \mathcal{C}_p$ is represented by a strictly positive continuous function.  Define:
\begin{align*}
  \widetilde{\vol}(f): =\inf_{K\in \mathfrak{K}^n} \left(\frac{f \cdot S(K^{n-1};\cdot)}{\vol(K)^{n-1/n}}\right)^{n}.
\end{align*}
where $f\cdot S(K^{n-1};\cdot)$ is the paring of two dual vector spaces.  This is a homogenous version of the classical polar transform in convex geometry.  Note that by the Brunn-Minkowski inequality, if $f=h_L$ for some $L\in \mathfrak{K}^n$ then $\widetilde{\vol}(h_L)=\vol(L)$.

We claim that this definition coincides with the definition using the Alexandrov decomposition described in Section \hyperref[section alex body]{4}.  Indeed, let $f=P(f)+N(f)$ be the decomposition of a strictly positive continuous function.  It is clear that $$\widetilde{\vol}(f)\ge \widetilde{\vol}(P(f))=\vol(P(f)).$$
Using the orthogonality of the decomposition, we also have $$\widetilde{\vol}(f) \leq  \left(\frac{f \cdot \langle P(f)^{n-1} \rangle}{\vol(P(f))^{n-1/n}}\right)^{n} = \vol(P(f)).$$  This establishes the equality.

\section{Blaschke addition for (1,1)-classes and its applications}\label{section blachke}
In this section, inspired by Blaschke addition for convex bodies, we introduce Blaschke addition for (1,1)-classes and give some applications.  Throughout we work with projective varieties.

We start with the convex geometry setting (see \cite[Chapter 8]{schneiderconvex}).  For convex bodies $K, L \in \mathfrak{K}^n$, Minkowski's existence theorem yields a convex body $M\in \mathfrak{K}^n$, unique up to translation, such that
\begin{equation*}
  S(M^{n-1};\cdot)=S(K^{n-1};\cdot)+S(L^{n-1};\cdot).
\end{equation*}

By the discussions in the previous sections, movable curve classes correspond to mixed area measures. Let $X$ be a smooth projective variety of dimension $n$. Assume that $L, M$ are two big (1,1)-classes.  As proved in \cite{lehmann2016positiivty}, there exists a unique big and movable (1,1)-class $N$ such that
\begin{equation*}
  \langle L^{n-1}\rangle +\langle M^{n-1}\rangle
  =\langle N^{n-1}\rangle.
\end{equation*}

Now we give the following definition.

\begin{defn}[\textsf{Blaschke addition}] \rm
Let $L, M$ be two big (1,1)-classes.  The Blaschke addition of $L, M$, denoted by $L \# M$, is defined to be the unique big and movable (1,1)-class satisfying
\begin{equation*}
\langle (L \# M)^{n-1}\rangle = \langle L^{n-1}\rangle +\langle M^{n-1}\rangle.
\end{equation*}
\end{defn}

\begin{rmk} \rm
By the definition of positive products, we know that $\langle L^{n-1}\rangle = \langle P(L)^{n-1}\rangle$.  Thus Blaschke addition of arbitrary big classes only depends on their positive parts.  For example, for any big divisor class $L$ we have $L \# L = 2^{\frac{1}{n-1}}P(L)$.
\end{rmk}

\begin{rmk} \rm
As pointed out to us by M.~Jonsson, the Blaschke addition can also defined on Riemann-Zariski spaces. More precisely, we can define the Blaschke addition of two Cartier big $\mathbf{b}$-divisor classes as a Weil movable $\mathbf{b}$-divisor class.
Let $X$ be a smooth projective variety of dimension $n$ and let $\mathfrak{X}$ be the Riemann-Zariski space containing all the smooth birational models of $X$. Assume that $\mathbf{L}=\{L_\alpha\}, \mathbf{M}=\{M_\alpha\}$ are big Cartier $\mathbf{b}$-divisor classes on $\mathfrak{X}$. Let $\varphi_{\alpha \beta}: X_\alpha \rightarrow X_\beta$ be a birational morphism.  We claim that $\varphi_{\alpha \beta *} (L_\alpha \# M_\alpha) = L_\beta \# M_\beta$, thus $\mathbf{L} \# \mathbf{M}$ is a Weil movable $\mathbf{b}$-divisor class. To this end, note that by the proof of Proposition \ref{pullback} and the definition of Cartier divisor class we have
\begin{align*}
\varphi_{\alpha \beta} ^* (\langle L_\beta ^{n-1}\rangle + \langle M_\beta ^{n-1}\rangle) &= \langle \varphi_{\alpha \beta} ^*L_\beta ^{n-1}\rangle + \langle \varphi_{\alpha \beta} ^*M_\beta ^{n-1}\rangle\\
&=\langle L_\alpha ^{n-1}\rangle + \langle M_\alpha ^{n-1}\rangle\\
&=\langle (L_\alpha \# M_\alpha) ^{n-1}\rangle,
\end{align*}
which implies that $\langle (L_\alpha \# M_\alpha) ^{n-1}\rangle = \langle \varphi_{\alpha \beta} ^*(L_\beta \# M_\beta) ^{n-1}\rangle$. Then we get $L_\alpha \# M_\alpha = \langle \varphi_{\alpha \beta} ^*(L_\beta \# M_\beta) \rangle$.  Taking $\varphi_{\alpha \beta *}$ of both sides yields $\varphi_{\alpha \beta *}(L_\alpha \# M_\alpha)= \langle L_\beta \# M_\beta \rangle = L_\beta \# M_\beta$ (as $L_\beta \# M_\beta$ is movable).
\end{rmk}

Following \cite{firey67mixedbody, lutwak1986mixedbody},
we can define the mixed divisor class associated to $(n-1)$ big $(1,1)$-classes.

\begin{defn}[\textsf{mixed divisor class}] \rm
Let $L_1, ..., L_{n-1}$ be big (1,1)-classes on $X$.  The mixed (1,1)-class, denoted by $[L_1, ..., L_{n-1}]$, is defined to be the unique big and movable (1,1)-class satisfying
\begin{equation*}
  \langle [L_1, ..., L_{n-1}]^{n-1} \rangle = \langle L_1 \cdot ... \cdot L_{n-1}\rangle.
\end{equation*}
\end{defn}

Before going on, we first give the following generalized Teissier proportionality theorem for big and movable $(1,1)$-classes.

\begin{thrm}\label{thm teisprop}
Let $X$ be a smooth projective variety of dimension $n$. Assume $L_1,..., L_n$ are big and movable $(1,1)$-classes on $X$, then \begin{align*}
  \langle L_1 \cdot ... \cdot L_n\rangle \geq \vol(L_1)^{1/n}...\vol(L_n)^{1/n}
\end{align*}
where the equality holds if and only if all the $L_i$ are proportional.
\end{thrm}

\begin{proof}
Firstly, by taking suitable Fujita approximations, it is easy to see that the inequality follows easily from the usual Khovanski-Teissier inequality for nef classes, hence we only need to characterize the equality situation.
Without loss of generality, we can assume that $\vol(L_i)=1$ for all $L_i$, then we need to show that
\begin{align}\label{eq teis eq}
  \langle L_1 \cdot ... \cdot L_n\rangle =1
\end{align}
if and only if $L_1 =L_2=...=L_n$. To this end, we claim that the equality (\ref{eq teis eq}) implies $\langle L_1 ^2 \cdot L_3 ... \cdot L_n\rangle =1$. Then by induction, we get that
\begin{align*}
  \langle L_1 ^{n-1} \cdot L_n\rangle =1.
\end{align*}
By the Teissier proportionality theorem proved in \cite{lehmann2016positiivty}, we conclude that $L_1 =L_n$.

The proof of our claim follows from the Hodge-Riemann bilinear relations. Consider the following quadratic form $Q$ on $H^{1,1}(X, \mathbb{R})$:
$$Q(\lambda, \mu):=  \lambda \cdot \mu \cdot A_1\cdot...\cdot A_{n-2},$$
where the $A_i$ are K\"ahler classes.
According to the Hodge-Riemann bilinear relations (see e.g. \cite[Theorem A and Theorem C]{DN06}), $Q$ is of signature $(1, \dim H^{1,1}(X, \mathbb{R}) -1)$. For any two nef classes $\alpha, \beta$, by considering the discriminant of the quadratic polynomial $p(t):= Q(\alpha+t\beta, \alpha+t\beta)$ and the signature of $Q$, we get
\begin{align}
\label{eq KT mixed HodgeRiemann}
(\alpha \cdot \beta \cdot A_1\cdot ...\cdot A_{n-2})^2
\geq ( \alpha^2 \cdot A_1\cdot ...\cdot A_{n-2})
\cdot (\beta^2 \cdot A_1\cdot ...\cdot A_{n-2}).
\end{align}
Since $\alpha, \beta, A_1, ..., A_{n-2}$ are arbitrary, we find by taking suitable Fujita approximations that
\begin{align*}
  \langle L_1 \cdot ... \cdot L_n\rangle ^2 \geq \langle L_1 ^2 \cdot L_3 ... \cdot L_n\rangle \cdot \langle L_2 ^2 \cdot L_3 ... \cdot L_n\rangle.
\end{align*}

Note that, in the above inequality, $\langle L_1 \cdot ... \cdot L_n\rangle =1$ and $\langle L_j ^2 \cdot L_3 ... \cdot L_n\rangle \geq 1$ for $j=1,2$. Hence we must have equality everywhere, in particular, we must have $\langle L_1 ^2 \cdot L_3 ... \cdot L_n\rangle =1$. This then finishes the proof of our claim.
\hfill $\Box$
\end{proof}

\begin{rmk} \rm
The above proof provides an alternative proof of a related result in \cite{lehmxiao16convexity}, where we applied complex Monge-Amp\`ere equations in big classes and pointwise Brunn-Minkowski inequalities for positive $(1,1)$-forms.
\end{rmk}

For Blaschke addition of big (1,1)-classes and mixed (1,1)-classes we obtain geometric inequalities similar to \cite{lutwak1986mixedbody}.

\begin{prop}\label{mixed class volbd}
Let $X$ be a smooth projective variety of dimension $n$. Assume $L_1,..., L_{n-1}$ are big and movable $(1,1)$-classes on $X$.  Then the volume of the mixed (1,1)-class $[L_1, ..., L_{n-1}]$ satisfies
\begin{equation*}
  \vol([L_1, ..., L_{n-1}])^{n-1} \geq \vol(L_1)\cdot...\cdot \vol(L_{n-1})
\end{equation*}
where the equality holds if and only if $L_{1},\ldots,L_{n-1}$ are proportional.
\end{prop}

\begin{proof}
This is a direct consequence of the Brunn-Minkowski inequality for (1,1)-classes. Note that we have $\vol([L_1, ..., L_{n-1}]) = \langle [L_1, ..., L_{n-1}]^{n-1} \rangle \cdot [L_1, ..., L_{n-1}]$, then
\begin{align*}
  \vol([L_1, ..., L_{n-1}])&= \langle L_1 \cdot ... \cdot L_{n-1}\rangle \cdot [L_1, ..., L_{n-1}] \\
  & \geq \langle L_1 \cdot ... \cdot L_{n-1}\cdot [L_1, ..., L_{n-1}]\rangle\\
  & \geq \vol(L_1)^{1/n} \cdot ...\cdot \vol(L_{n-1})^{1/n} \cdot \vol([L_1, ..., L_{n-1}])^{1/n},
\end{align*}
which implies the desired inequality.  The characterization of the equality follows from Theorem \ref{thm teisprop}.~\hfill $\Box$
\end{proof}

We use $[L, M]_i$ to denote the mixed (1,1)-class of
$\langle L^{n-1-i} \cdot M^i \rangle$. It is clear that $[L, M]_0 =L$ and $[L, M]_{n-1} =M$.

\begin{prop}
Let $X$ be a smooth projective variety of dimension $n$. Assume that $L, M$ are big and nef (1,1)-classes.  Then the sequence $\{\vol([L,M]_i)\}_{i=0} ^{i=n-1}$ is log-concave, that is,
\begin{equation*}
  \vol([L,M]_k)^2 \geq \vol([L,M]_{k-1}) \vol([L,M]_{k+1})\ \ \textrm{for any}\ k=1,...,n-2.
\end{equation*}
Moreover, the following statements are equivalent:
\begin{enumerate}
  \item[\rm 1.] $\vol([L,M]_k)^2 = \vol([L,M]_{k-1}) \vol([L,M]_{k+1})\ \ \textrm{for all}\ k=1,...,n-2$;
  \item[\rm 2.] $\vol([L,M]_{n-2})^{n-1} = \vol([L,M]_{0}) \vol([L,M]_{n-1})^{n-2}$;
  \item[\rm 3.] $L,M$ are proportional.
\end{enumerate}
\end{prop}

\begin{proof}
By the Hodge-Riemann bilinear relations as in the proof of Theorem \ref{thm teisprop},  for any movable (1,1)-class $N$ we find that
\begin{equation*}
  (L\cdot M\cdot L^{n-2-k}\cdot M^{k-1} \cdot N)^2 \geq ( L^2 \cdot L^{n-2-k}\cdot M^{k-1}\cdot N)( M^2\cdot L^{n-2-k}\cdot M^{k-1}\cdot N ).
\end{equation*}

By the definition of mixed class and Brunn-Minkowski inequalities, if we take $N= [L, M]_k$, then we get that
\begin{align*}
  \vol([L,M]_k)^2 &\geq (\langle [L, M]_{k-1} ^{n-1}\rangle \cdot [L,M]_k) (\langle [L, M]_{k+1} ^{n-1}\rangle \cdot [L,M]_k) \\
  & \geq \vol([L, M]_{k-1})^{n-1/n}\vol([L, M]_{k+1})^{n-1/n} \vol([L,M]_k)^{2/n}.
\end{align*}
This implies the log concavity: $\vol([L,M]_k)^2 \geq \vol([L,M]_{k-1}) \vol([L,M]_{k+1})$.

For the equality characterizations, first note that the equivalence of $(2)$ and $(3)$ is a direct consequence of Proposition \ref{mixed class volbd}. We only need to verify the equivalence of $(1)$ and $(2)$, but this follows directly from the log concavity of the sequence $\{\vol([L,M]_i)\}_{i=0} ^{i=n-1}$.
\hfill $\Box$
\end{proof}

By the properties of positive product, Blaschke addition is compatible with mixed divisor class construction in the following sense.

\begin{prop}\label{BlasAddprop}
Let $X$ be a projective variety of dimension $n$.  Assume that $L, M$ are big and nef. For any big (1,1)-classes $D_2, ..., D_{n-1}$ we have
\begin{equation*}
  [L+M, D_2, ..., D_{n-1}]= [L, D_2, ..., D_{n-1}]\# [M, D_2, ..., D_{n-1}]
\end{equation*}
where $\#$ is the Blaschke addition.
\end{prop}

\begin{proof}
Firstly, by the definition of positive products, we always have \begin{align*}
  \langle (L+M)\cdot D_2 \cdot...\cdot D_{n-1} \rangle \succeq \langle L\cdot D_2 \cdot...\cdot D_{n-1} \rangle + \langle M\cdot D_2 \cdot...\cdot D_{n-1} \rangle.
\end{align*}
On the other hand, since $L, M$ are nef, we claim that the equality holds: for any ample class $A$, we have
\begin{align*}
  \langle (L+M)\cdot D_2 \cdot...\cdot D_{n-1} \rangle\cdot A
  &= \langle (L+M)\cdot D_2 \cdot...\cdot D_{n-1} \cdot A \rangle\\
  &=  (L+M)\cdot\langle D_2 \cdot...\cdot D_{n-1} \cdot A \rangle\\
  &= \langle L\cdot D_2 \cdot...\cdot D_{n-1} \rangle \cdot A + \langle M\cdot D_2 \cdot...\cdot D_{n-1} \rangle \cdot A.
\end{align*}

By the definition of mixed class, we have
\begin{align*}
  \langle ([L, D_2, ..., D_{n-1}]\# [M, D_2, ..., D_{n-1}])^{n-1}\rangle &= \langle [L, D_2, ..., D_{n-1}]^{n-1}\rangle + \langle [M, D_2, ..., D_{n-1}]^{n-1}\rangle\\
  &= \langle L\cdot D_2 \cdot...\cdot D_{n-1} \rangle + \langle M\cdot D_2 \cdot...\cdot D_{n-1} \rangle \\
  &=\langle (L+M)\cdot D_2 \cdot...\cdot D_{n-1} \rangle\\
  &=\langle [L+M, D_2, ..., D_{n-1}]^{n-1}\rangle.
\end{align*}

As a consequence of the Teissier proportionality for big and movable classes proved in \cite{lehmann2016positiivty}, we immediately get the equality $[L+M, D_2, ..., D_{n-1}]= [L, D_2, ..., D_{n-1}]\# [M, D_2, ..., D_{n-1}]$.
\hfill $\Box$
\end{proof}

Similar to the convex geometry setting \cite{lutwak1986mixedbody}, using Blaschke addition and mixed divisor classes, we get a somewhat improved Khovanskii-Teissier inequality for (1,1)-classes.

\begin{thrm}\label{thrm improved BM}
Let $X$ be a smooth projective variety of dimension $n$. Assume that $L, M$ are big and nef (1,1)-classes, then
\begin{align*}
  \vol(L+M)^{1/n} &\geq \left(\sum_{i=0} ^{n-1} \frac{(n-1)!}{i!(n-1-i)!} \vol([L, M]_i)^{n-1/n}\right)^{1/n-1}\\
  &\geq \vol(L)^{1/n} + \vol(M)^{1/n}.
\end{align*}
The equality in either inequality holds if and only if $L, M$ are proportional.
\end{thrm}

\begin{proof}
By Proposition \ref{BlasAddprop} we have
\begin{equation}\label{eq mixedbi}
  [L+M,..., L+M] = \sum_{i=0} ^{n-1} \frac{(n-1)!}{i!(n-1-i)!}\cdot [L, M]_i,
\end{equation}
where the above addition is the Blaschke addition.
We claim that
\begin{equation} \label{ineq mixedclass}
  \vol(L+M)^{n-1/n} \geq \sum_{i=0} ^{n-1} \frac{(n-1)!}{i!(n-1-i)!} \vol([L, M]_i)^{n-1/n}
\end{equation}
with equality if and only if $L, M$ are proportional. To this end, note that
\begin{align*}
  \vol(L+M) &= \langle [L+M,...,L+M]^{n-1}\rangle\cdot (L+M) \\
  & = \langle (\sum_{i=0} ^{n-1} \frac{(n-1)!}{i!(n-1-i)!}\cdot [L, M]_i ) ^{n-1}\rangle \cdot (L+M)\\
  & =\sum_{i=0} ^{n-1} \frac{(n-1)!}{i!(n-1-i)!} \langle [L, M]_i  ^{n-1}\rangle \cdot (L+M)
\end{align*}
where the third line follows from the definition of Blaschke addition. Then by the Brunn-Minkowski inequality, we immediately get the claimed inequality. If the equality holds, then we must have
\begin{equation*}
  \langle [L, M]_i  ^{n-1}\rangle \cdot (L+M) = \vol([L, M]_i)^{n-1/n}\vol(L+M)^{1/n}\ \ \textrm{for every}\ i.
\end{equation*}
In particular, applying Theorem \ref{thm teisprop} to the cases when $i=0, n-1$, we get that $L, M$ are proportional.

On the other hand, by Proposition \ref{mixed class volbd} we have
\begin{equation*}
  \vol([L, M]_i)^{n-1} \geq \vol(L)^{n-1-i} \vol(M)^i .
\end{equation*}
This implies that
\begin{equation}\label{ineq2 mixedclass}
  \sum_{i=0} ^{n-1} \frac{(n-1)!}{i!(n-1-i)!} \vol([L, M]_i)^{n-1/n} \geq \left(\vol(L)^{1/n} + \vol(M)^{1/n}\right)^{n-1}.
\end{equation}
For the characterization of the equality, note that the equality implies that, for every $i$ we have
\begin{equation*}
  \vol([L, M]_i)^{n-1} = \vol(L)^{n-1-i} \vol(M)^i.
\end{equation*}
Hence, by Proposition \ref{mixed class volbd}, we get that $L, M$ must be proportional.
\hfill $\Box$
\end{proof}

\subsection{Kneser-S\"{u}ss inequality for divisors}
In the convex geometry setting, let $K, L\in \mathfrak{K}^n$ be convex bodies, and let $M$ be the convex body satifying
\begin{equation*}
  S(M^{n-1};\cdot)=S(K^{n-1};\cdot)+S(L^{n-1};\cdot).
\end{equation*}
It is proved in \cite{kneser1932volumina} (see also \cite[Chapter 8]{schneiderconvex}) that
\begin{equation*}
  \vol(M)^{n-1/n} \geq \vol(K)^{n-1/n} +\vol(L)^{n-1/n},
\end{equation*}
where equality holds if and if $K, L$ are homothetic. This inequality is called the Kneser-S\"{u}ss inequality.

Let $X$ be a smooth projective variety of dimension $n$, in \cite{xiao15} we defined the following volume type function for movable curve classes:
\begin{align*}
  \mathfrak{M}(\alpha)=\inf_{B\ \textrm{big divisor}} \left(\frac{\alpha\cdot B}{\vol(B)^{1/n}} \right)^{n/n-1},\ \alpha\in \Mov_1 (X).
\end{align*}
By its definition, we have
\begin{equation}\label{eq kneser}
 \mathfrak{M}(\alpha_1 +\alpha_2)^{n-1/n} \geq  \mathfrak{M}(\alpha_1)^{n-1/n} + \mathfrak{M}(\alpha_2)^{n-1/n}.
\end{equation}
Later, in \cite{lehmann2016positiivty} (see also \cite{lehmxiao16convexity}), by analyzing ``Zariski decomposition structure'' of $\mathfrak{M}$, we have

\begin{thrm}
Let $X$ be a smooth projective variety of dimension $n$, and assume that $\alpha_i \in \Mov_1 (X)$ satisfies $\mathfrak{M}(\alpha_i)>0$ ($i=1, 2$), then
\begin{equation}
 \mathfrak{M}(\alpha_1 +\alpha_2)^{n-1/n} \geq  \mathfrak{M}(\alpha_1)^{n-1/n} + \mathfrak{M}(\alpha_2)^{n-1/n}
\end{equation}
with equality if and only if $\alpha_1 ,\alpha_2$ are proportional.
\end{thrm}

Equivalently, if we translate the above statement to Blaschke addition for $(1,1)$-classes, then we get
\begin{thrm}
Let $X$ be a smooth projective variety of dimension $n$, and assume that $L, M$ are big movable $(1,1)$-classes, then
\begin{equation}
 \vol(L\#M)^{n-1/n} \geq \vol(L)^{n-1/n} +\vol(M)^{n-1/n}
\end{equation}
with equality if and only if $L, M$ are proportional.
\end{thrm}

\section{Projection volumes and the Brascamp-Lieb inequality}
\label{section loomis}
There are many geometric inequalities relating the volume of a convex body to the volumes of its images under linear maps.  One of the most well-known is the Loomis-Whitney inequality \cite{loomis1949inequality}: let $K$ be a measurable subset of $\mathbb{R}^n$, and let $\pi_j : \mathbb{R}^n \rightarrow \mathbb{R}^{n-1},\ (x_1,..., x_n) \mapsto \widehat{x}_j:=(x_1, ..., x_{j-1}, x_{j+1},..., x_n)$ be projections, then
\begin{equation*}
  \vol(K)^{n-1} \leq \prod_{j=1} ^{n} \vol(\pi_j (K)).
\end{equation*}
More generally, the Loomis-Whitney inequality finds its natural generalisation in the Brascamp-Lieb inequality for general linear projections (see e.g. \cite{braslieb76}, \cite{tao2008brascamp}), but the inequality needs to be corrected by a factor determined by the projections.

In this section, we formulate and prove an analogue for $(1,1)$ classes.

\subsection{Projection volumes}

Fix a convex body $P \subseteq \mathbb{R}^{n}$ and a unit vector $v \in \mathbb{R}^{n}$.  Then the volume of the orthogonal projection of $P$ onto $v^{\perp}$ is computed by the mixed volume $n \langle P^{n-1},  \overline{0v} \rangle$, where $\overline{0v}$ is the segment between $0$ and $v$.  The volume of the projection onto the orthogonal of a subspace $W$ can be computed by a similar procedure.

The analogous construction in complex geometry is as follows.  In place of an orthogonal projection, we fix a dominant rational map $f: X \dashrightarrow Z$ whose resolution has connected fibers.  (This is motivated by toric geometry, where a ``direction'' in the polytope space identifies a map on a blow-up.)  Recall that \cite{BFJ09} shows that for a big divisor $L$ and a prime divisor $H$ the restricted volume $\vol_{X|H} (L)$ is just the positive product $\langle L^{n-1}\rangle\cdot H$. Given a big class $L$, the projection volume onto the fiber is defined to be
\begin{equation*}
\vol_{f}(L) := \langle \phi^{*}L^{n-d} \rangle \cdot F
\end{equation*}
where $d=\dim Z$, $\phi: X' \to X$ is a resolution of $f$ and $F$ is a general fiber of the induced map on $X'$.
We next need to ask when a collection $\{f_{i}: X \dashrightarrow Z_{i}\}$ is ``mutually orthogonal''.  Guided by the toric case, we require that the induced rational map
$X \dashrightarrow \prod_{i} Z_{i}$
is birational.

\subsection{Loomis-Whitney inequality}
Inspired by the Loomis-Whitney inequality for convex bodies, we first prove the following analogue on projective varieties.

\begin{thrm}[Loomis-Whitney inequality]\label{thrm loomis}
Let $X$ be a smooth projective variety of dimension $n$.  Suppose that $\{ f_{i}: X \dashrightarrow C_{i} \}_{i=1}^{n}$ is a collection of dominant rational maps onto curves whose resolution has connected fibers, and assume that the induced map $X \dashrightarrow \prod_{i=1} ^n C_{i}$ is birational. Then for any big $(1,1)$ class $L$, we have
\begin{equation*}
\left( \vol(L)/n! \right)^{n-1}  \leq \prod_{i=1}^{n} \left( \vol_{f_i} (L)/(n-1)! \right).
\end{equation*}
If all the $f_{i}$ are morphisms, then we obtain equality if and only if the positive part of $L$ is numerically proportional to a sum of fibers of the $f_{i}$.
\end{thrm}

Note that the inclusion of the factorials corrects for the discrepancy between volumes in the complex geometry and in the convex geometry settings.

The first proof of the Loomis-Whitney inequality in the convex geometry setting was given by \cite{loomis1949inequality}.  The strategy is to approximate a polytope by hypercubes and to use a counting argument which inducts on dimension. In our setting, we follow the approach of \cite{leng07loomismixed}, where the authors mainly applied the method of \cite{zhang99sobolev} using mixed volumes.

\begin{proof}[Proof of Theorem \ref{thrm loomis}]\label{rankone BL}
After passing to a birational model of $X$ and pulling back $L$, we may suppose that each rational map $f_{i}$ is a morphism.

Denote the general fiber of $f_i$ by $F_i$, and assume that $a_1, ..., a_n$ are positive numbers. Note

\begin{align*}
 \vol(\sum_{i=1} ^n a_i F_i) =  n! \prod_{i=1} ^n a_i.
\end{align*}

By the Khovanski-Teissier inequality for big divisors  as in \cite{lehmann2016positiivty}, we get
\begin{align*}
\langle L^{n-1}\rangle \cdot (\sum_{i=1} ^n a_i F_i) &\geq \vol(L)^{n-1/n} \vol(\sum_{i=1} ^n a_i F_i)^{1/n}\\
& \geq  (n!\prod_{i=1} ^n a_i)^{1/n}\vol(L)^{n-1/n}
\end{align*}
with equality iff the positive parts of the divisors are numerically proportional.
Setting $a_i = c_i / \langle L^{n-1}\rangle \cdot F_i$ for a variable $c_{i}$, the above inequality implies that
\begin{equation*}
  \prod_{i=1} ^n \langle L^{n-1}\rangle \cdot F_i \geq n! \vol(L)^{n-1} \frac{\prod_{i=1} ^n c_i}{(\sum_{i=1} ^n c_i)^n}.
\end{equation*}
Note that $\frac{\prod_{i=1} ^n c_i}{(\sum_{i=1} ^n c_i)^n} \leq 1/n^n$ with equality if and only if all the $c_i$ are equal. Letting $c_i =1$, we get
\begin{equation*}
  \prod_{i=1} ^n \langle L^{n-1}\rangle \cdot F_i /(n-1)! \geq (\vol(L)/n!)^{n-1}.
\end{equation*}
Furthermore, equality holds here if and only if equality holds in the Khovanskii-Teissier inequality above, finishing the proof of our theorem.
\hfill $\Box$
\end{proof}

\subsection{Rank-one Brascamp-Lieb inequality}

The following statement can be seen as a special case of the Brascamp-Lieb inequality.  In the convex geometry setting, it corresponds to comparing the volume of a polytope to its projections onto the coordinate axes (instead of the coordinate hyperplanes).

\begin{prop} [Rank-one Brascamp-Lieb inequality]
\label{rankoneBL}
Suppose that $X$ is a smooth projective variety birational to a product of curves $\prod_{i=1}^{n} C_{i}$.  Then for any big class $L$ on $X$ we have
\begin{equation*}
\vol(L)/n! \leq \prod_{i=1}^{n} \langle L \rangle \cdot  T_{i}
\end{equation*}
where $T_{i}$ is a general fiber of the rational map $X \dashrightarrow \prod_{j \neq i} C_{j}$.
\end{prop}

\begin{proof}
The proof is by induction on the dimension $n$.

By passing to a birational model of $X$ and pulling back $L$, we may suppose that each rational map $f_{i}$ is a morphism. Note that the quantities $\langle L \rangle \cdot T_{i}$ do not change upon such pullback, since a very general fiber $T_{i}$ avoids the restricted base locus of $\langle L \rangle$.  We may also replace $L$ by $\langle L \rangle$ so that we may assume $L$ is big and movable.

Let $F_{i}$ denote a very general fiber of the map $X \to C_{i}$.  Let $s$ denote the largest real number such that $L - sF_{1}$ is pseudo-effective.  Note that $s \leq L \cdot T_{1}$ since $F_{1} \cdot T_{1} = 1$ and $(L - sF_{1}) \cdot T_{1} \geq 0$.  Then we have
\begin{equation*}
\vol(L)/n = \int_{0}^{s} \langle (L - tF_{1})^{n-1} \rangle \cdot F_{1} \, dt \leq (L \cdot T_{1}) (\langle L^{n-1} \rangle \cdot F_{1}).
\end{equation*}
Clearly $\langle L^{n-1} \rangle \cdot F_{1}$ is at most $\vol(L|_{F_{1}})$.  And by induction, we have
\begin{equation*}
 \vol(L|_{F_{1}})/(n-1)! \leq \prod_{i=2}^{n} \langle L|_{F_{1}} \rangle \cdot T_i.
\end{equation*}
Furthermore, for $2 \leq i \leq n$ we have $\langle L|_{F_{1}} \rangle \cdot T_{i} = L \cdot T_{i}$.  Indeed, this follows from the fact that $L|_{F_{1}}$ is movable when $F_{1}$ is very general, so that the positive product on the left agrees with $L|_{F_{1}}$.  We conclude that
\begin{equation*}
\vol(L)/n \leq (n-1)! \prod_{i=1}^{n} L \cdot T_{i}
\end{equation*}
giving the statement.
\hfill $\Box$
\end{proof}

\providecommand{\bysame}{\leavevmode\hbox to3em{\hrulefill}\thinspace}
%
%

\bibliographystyle{amsalpha}
\bibliographymark{References}

\end{document}